\documentclass[acmtog, authorversion=true]{acmart}

\copyrightyear{2024}
\acmYear{2024}
\setcopyright{acmlicensed}\acmConference[SA Conference Papers '24]{SIGGRAPH Asia 2024 Conference Papers}{December 3--6, 2024}{Tokyo, Japan}
\acmBooktitle{SIGGRAPH Asia 2024 Conference Papers (SA Conference Papers '24), December 3--6, 2024, Tokyo, Japan}
\acmDOI{10.1145/3680528.3687599}
\acmISBN{979-8-4007-1131-2/24/12}

\usepackage{graphicx} %
\usepackage{cleveref}
\usepackage{subcaption}
\usepackage{rotating}
\usepackage[ruled]{algorithm2e}
\usepackage{hyperref}

\crefformat{section}{Section~#2#1#3}
\crefformat{subsection}{Section~#2#1#3}
\crefformat{equation}{Eq.~#2#1#3}
\crefformat{figure}{Fig.~#2#1#3}
\crefformat{algocf}{Alg.~#2#1#3}

\newcommand{\dimension}{d}
\newcommand{\sol}{u}
\newcommand{\Rd}{{\mathbb{R}^\dimension}}
\newcommand{\Rtwo}{{\mathbb{R}^2}}
\newcommand{\Rthree}{{\mathbb{R}^3}}

\newcommand{\Surface}{\mathcal{S}}
\newcommand{\Boundary}{\mathcal{C}}
\newcommand{\Medial}{\mathrm{med}}
\newcommand{\Tube}{\mathcal{N}}
\newcommand{\lfs}{\mathrm{LFS}}
\newcommand{\mylaplace}{\Delta}
\newcommand{\mylaplacebel}{\mylaplace_\Surface}
\newcommand{\grad}{\nabla}

\renewcommand{\vec}{\mathbf}
\newcommand{\vecx}{\vec{x}}
\newcommand{\vecy}{\vec{y}}
\newcommand{\vecz}{\vec{z}}
\newcommand{\vecr}{\vec{r}}
\newcommand{\vecq}{\vec{q}}
\newcommand{\vecn}{\vec{n}}
\newcommand{\vect}{\vec{t}}
\newcommand{\vech}{\vec{h}}
\newcommand{\radius}{r}

\newcommand{\dAy}{\,\mathrm{d}\vecy}

\newcommand{\dVz}{\,\mathrm{d}\vecz}
\newcommand{\source}{f}
\newcommand{\cps}{\mathrm{cp}_\Surface}
\newcommand{\cpc}{\mathrm{cp}_\Boundary}
\newcommand{\estimate}{\widehat}
\newcommand{\distance}{\delta}

\definecolor{S-blue}{rgb}{0.161,0.404,0.690}
\definecolor{grey}{rgb}{0.529, 0.529, 0.529}
\definecolor{dark-grey}{rgb}{0.3, 0.3, 0.3}

\defcitealias{RiouxLavoie2022:Fluids}{Rioux-Lavoie and Sugimoto et~al.}
\defcitealias{Sawhney:2022:GFMC}{Sawhney and Seyb et~al.}
\defcitealias{sawhney2023walk}{Sawhney and Miller et~al.}
\defcitealias{SabelfeldSimonov1994}{Sabelfeld and Simonov}
\defcitealias{Miller2023}{Miller and Sawhney et~al.}
\defcitealias{Miller:2024:WRW}{Miller and Sawhney et~al.}

\newcommand{\mycite}[1]{[\citetalias{#1}~\citeyear{#1}]}
\newcommand{\mycitet}[1]{\citetalias{#1}~[\citeyear{#1}]}

\title{Projected Walk on Spheres: A Monte Carlo Closest Point Method for Surface PDEs}

\author{Ryusuke Sugimoto}
\orcid{0000-0001-5894-0423}
\affiliation{
    \institution{University of Waterloo}
    \country{Canada}
}
\email{rsugimot@uwaterloo.ca}
\author{Nathan King}
\orcid{0000-0003-4105-0189}
\affiliation{
    \institution{University of Waterloo}
    \country{Canada}
}
\email{n5king@uwaterloo.ca}
\author{Toshiya Hachisuka}
\orcid{0000-0003-0284-3776}
\affiliation{
    \institution{University of Waterloo}
    \country{Canada}
}
\email{toshiya.hachisuka@uwaterloo.ca}
\author{Christopher Batty}
\orcid{0000-0003-3830-7772}
\affiliation{
    \institution{University of Waterloo}
    \country{Canada}
}
\email{christopher.batty@uwaterloo.ca}



\begin{document}

\begin{abstract}
We present \emph{projected walk on spheres} (PWoS), a novel pointwise and discretization-free Monte Carlo solver for surface PDEs with Dirichlet boundaries, as a generalization of the \emph{walk on spheres} method (WoS)~\cite{Sawhney:2020:MCGP, Muller:1956:Some}. We adapt the recursive relationship of WoS designed for PDEs in volumetric domains to a volumetric neighborhood around the surface, and at the end of each recursion step, we project the sample point on the sphere back to the surface. We motivate this simple modification to WoS with the theory of the closest point extension used in the closest point method. 
To define the valid volumetric neighborhood domain for PWoS, we develop strategies to estimate the local feature size of the surface and to compute the distance to the Dirichlet boundaries on the surface extended in their normal directions. We also design a mean value filtering method for PWoS to improve the method's efficiency when the surface is represented as a polygonal mesh or a point cloud. Finally, we study the convergence of PWoS and demonstrate its application to graphics tasks, including diffusion curves, geodesic distance computation, and wave propagation animation. We show that our method works with various types of surfaces, including a surface of mixed codimension.\sloppy
\end{abstract}

\maketitle

\begin{figure}
\centering
\includegraphics[width=\linewidth]{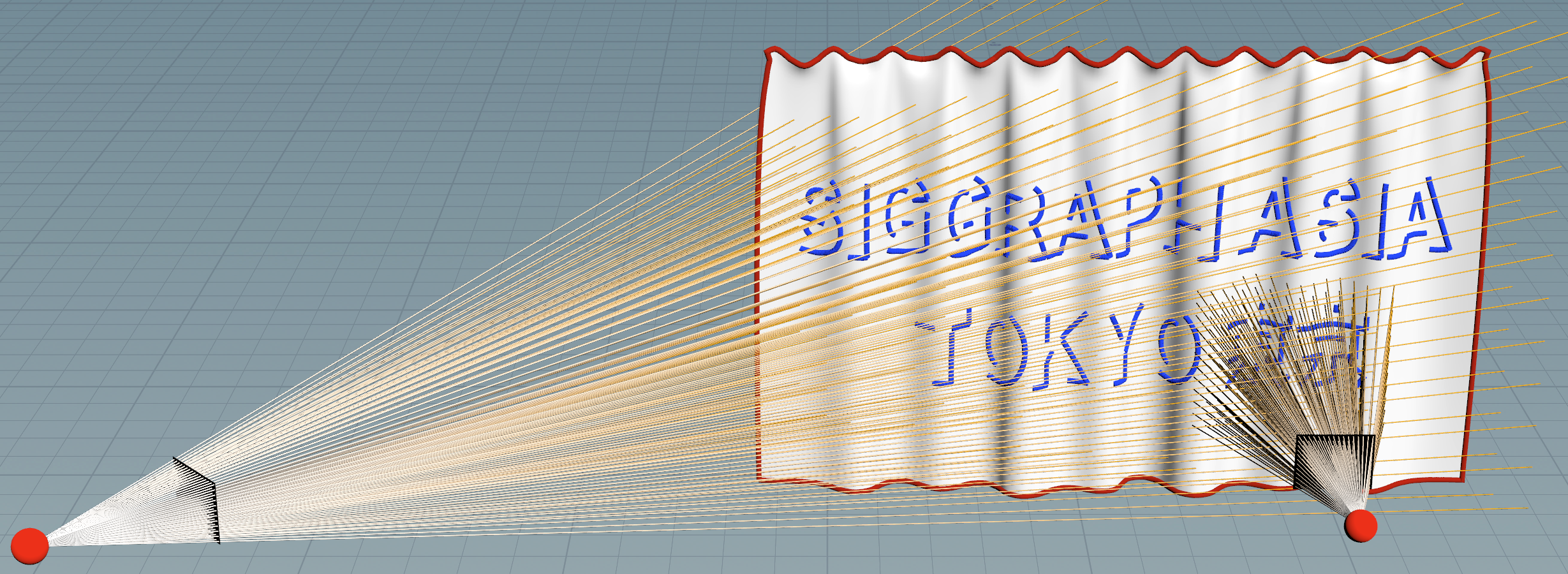}
\begin{minipage}[ht]{0.5\linewidth}\centering
{\setlength{\fboxsep}{0pt}\fbox{%
\includegraphics[width=\linewidth]{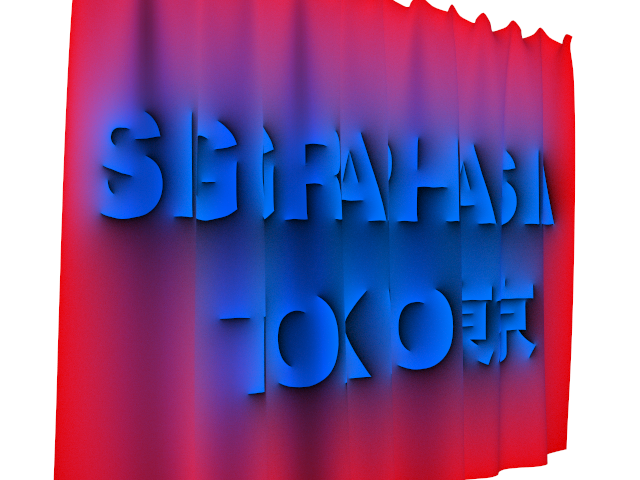}%
}}
\end{minipage}%
\begin{minipage}[ht]{0.5\linewidth}\centering
{\setlength{\fboxsep}{0pt}\fbox{%
\includegraphics[width=\linewidth]{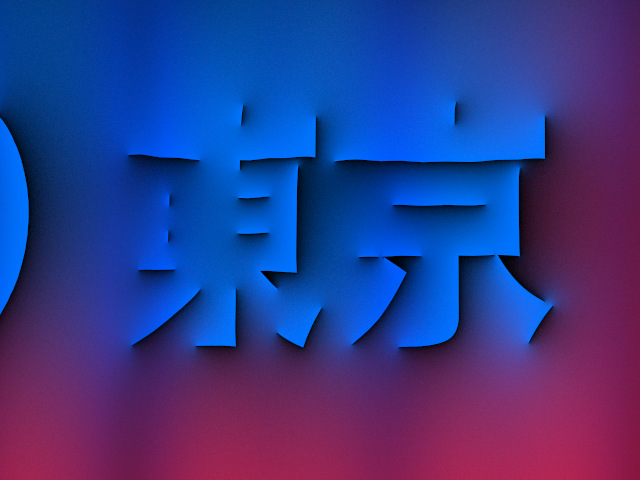}%
}}
\end{minipage}%
\caption{View-dependent diffusion curves with PWoS. Using our method, we solve the Laplace equation on a curved surface in a view-dependent manner. The pointwise and discretization-free nature of PWoS allows for the evaluation of the colors only at visible points where the object color is required by a shading algorithm with stochastic pixel-filtering.}
\label{fig:diffusioncurves_viewdependent}
\end{figure}

\section{Introduction}
Methods to solve \emph{surface partial differential equations} (PDEs) have become ubiquitous tools in computer graphics research and production. They are used for surface editing~\cite{Desbrun1999}, texture synthesis~\cite{Turk1991}, surface fluid animation~\cite{Stam2003}, geodesic distance computation~\cite{Crane2013}, and diffusion curves on surfaces~\cite{Bartsch2023, DeGoes2022} today. A common approach to tackle these problems is to discretize the surface and solve a globally coupled linear system using discrete surface differential operators.

For \emph{volumetric} PDEs, Monte Carlo methods have recently garnered significant attention in the graphics community due to their unique advantages over traditional discretization-based PDE solvers, including the ability to estimate solution values in a pointwise, spatial discretization-free manner. One such method is the \emph{walk on spheres} (WoS) method \cite{Muller:1956:Some} introduced to graphics by \citet{Sawhney:2020:MCGP}. They primarily focused on the (constant coefficient) Poisson and screened Poisson equations in a volumetric domain, and follow-up work likewise emphasizes volumetric problems. In the present work, we consider instead the problem of solving \emph{surface PDEs}.

\mycitet{Sawhney:2022:GFMC} proposed an extension of WoS to support second-order linear elliptic PDEs with spatially varying coefficients, and as one application, they demonstrated a method to solve the Laplace equation on a 2D surface embedded in 3D. However, this approach requires that the conformal parametrization of the surface be readily available, limiting the method's applicability. 

We propose a simpler generalization of WoS for surface PDEs, the \emph{projected walk on spheres} (PWoS) method, which only assumes the availability of a closest surface point query and an unoriented surface normal direction query. PWoS supports Dirichlet boundary conditions and inherits the advantages of WoS: PWoS is a pointwise, discretization-free Monte Carlo method.
Since our method does not require the meshing of the surface, it is particularly advantageous over traditional approaches, such as the finite element method, when the computation can be localized and when the surface is given as an implicit representation, such as a signed distance function.
The resulting solution is also free of mesh-dependent discretization artifacts, such as from linear interpolation,
as we show in \cref{fig:diffusioncurves_viewdependent}. 
Compared to WoS, which performs walks on spheres inside the domain, PWoS performs walks on spheres inside a Cartesian embedding neighborhood domain around the surface. After each step of the walk, it \emph{projects} the sampled point on the sphere onto the surface.
We motivate this simple modification to the original WoS through its connection to the \emph{closest point method} (CPM) \cite{Ruuth2008:CPM, Marz2012:CPM}.

Furthermore, inspired by the mean value filtering method for WoS by \citet{Bakbouk2023:MVC}, we design a mean value filtering method with a discrete basis function to allow more efficient estimation of solutions when the surface is represented as a polygonal mesh or a point cloud. To confirm PWoS's accuracy, we perform convergence studies of the method applied to the surface Poisson and screened Poisson equations. Finally, we demonstrate its use in several representative graphics applications, including diffusion curves, geodesic distance computation, and wave propagation animation.

In summary, our key contributions are:
\begin{itemize}
    \item Our novel PWoS algorithm that generalizes the WoS method to solve surface PDEs with Dirichlet boundaries, supported by the theory of the closest point extension.
    \item A mean value filtering method for PWoS with a discrete basis for efficiency improvement.
    \item Evaluation of PWoS with convergence studies and graphics applications.
\end{itemize}

\section{Background}
This section briefly reviews the two core mathematical ideas on which our method is based.

\subsection{Walk on Spheres}
WoS solves volumetric PDEs such as the Poisson equation over a Cartesian domain in $\Rd$.
Consider a $\dimension$-ball $B_\radius(\vecx)$, centered at $\vecx$ with radius $r$, fully contained within the domain. The integral equation 
\begin{equation}\label{eq:integraleq}
    \sol(\vecx) = \frac{1}{\lvert \partial B_\radius(\vecx) \rvert} \int_{\partial B_\radius(\vecx)} \sol(\vecy) \dAy + \int_{B_\radius(\vecx)} \source(\vecz) G(\vecx, \vecz) \dVz
\end{equation}
holds for the Poisson equation $\mylaplace u = \source$ in general, where $\lvert \partial B_\radius(\vecx) \rvert$ denotes the surface area of the sphere that bounds the ball $B_\radius(\vecx)$ and $G$ denotes the Green's function for the Poisson equation on $B_\radius(\vecx)$~\cite{Sawhney:2020:MCGP}. On the right-hand side, the first term is a boundary integral over the \mbox{$(\dimension-1)$-sphere}, and the second term is a volume integral over the $\dimension$-ball. 
If we perform Monte Carlo integration of the first term by uniformly sampling a point on the sphere and of the second term by sampling $N_V$ points $\vecz_i$ inside the ball with probability density function (PDF) $p(\vecz_i)$, we get the recursive relationship used in WoS:
\begin{equation}\label{eq:wos}
\estimate\sol(\vecx) = \estimate\sol(\vecy) + \frac{1}{N_V}\sum_{i=1}^{N_V}\frac{G(\vecx, \vecz_i) \source(\vecz_i)}{p(\vecz_i)},
\end{equation}
where the hat notation indicates that a variable is a Monte Carlo estimate. 
The first term on the right-hand side is also a Monte Carlo estimate because the solution $\sol(\vecy)$ is unknown at point $\vecy$ in general.
At each recursion step, WoS applies this recursive relationship to the largest ball inside the domain bounded by Dirichlet boundaries. It terminates the recursion when the sample point $\vecx$ during the recursion falls within a small distance $\epsilon$ of the domain boundary, by using the known solution at the closest boundary point $\overline\vecx$ as the solution estimate: \mbox{$\estimate\sol(\vecx) = \sol(\overline\vecx)$}. WoS thereby estimates the solution at each evaluation point independently, offering intrinsic parallelism. Our method generalizes WoS, originally proposed for volumetric PDEs, by incorporating the closest point extension theory of CPM.

\subsection{Surface PDEs and Closest Point Extension}\label{sec:cpm}
Consider the Poisson equation defined on a surface $\Surface$ embedded in $\Rd$ such that $\dim(\Surface) < \dimension$:
\begin{equation}\label{eq:surface_pde}
    \mylaplacebel \sol_\Surface (\vecy) = \source_\Surface (\vecy), \quad \vecy \in \Surface,
\end{equation}
where $\mylaplacebel$ denotes the Laplace-Beltrami operator. For convenience, we will use the word surface to refer to any nonzero codimension object in $\Rd$, including mixed-codimension objects. One typically solves such a surface PDE by discretizing the differential operator and solving a sparse linear system. For triangle meshes one can use the cotangent Laplacian~\cite{macneal1949solution}; for other surface representations, a corresponding discrete Laplacian must be defined~\cite{Bunge2023,Sharp2020,Belkin2009}. 
The closest point method~\cite{Ruuth2008:CPM} instead addresses the surface PDE (\cref{eq:surface_pde}) in a more general way by changing the domain to an embedding space, which is a Cartesian neighborhood surrounding the original surface. CPM then solves an \emph{embedding PDE} defined on the embedding space, whose solution when restricted to points $\vecy \in \Surface$ is the solution $u_\Surface(\vecy)$ to the original surface PDE. We briefly summarize CPM theory and refer readers to work by \citet{King:2023:CPM} for an in-depth review of CPM.

We first assume that $\Surface$ is smooth and define the closest point query to the surface, for $\vecx \in \Rd$, as
\begin{equation}\label{eq:cps}
    \cps(\vecx) = \underset{\vecy\in\Surface}{\operatorname{argmin}}\lVert \vecx - \vecy \rVert_2.
\end{equation}
In general, the mapping $\cps(\vecx)$ may not be unique: there may exist more than one closest point for a given $\vecx$. We define the neighborhood $\Tube(\Surface)$ where the closest point function is unique as
    $\Tube(\Surface) = \left\{\vecx \in \Rd \;\big|\; \lVert \vecx - \cps(\vecx) \rVert_2 < \lfs\left(\cps(\vecx)\right)\right\}$, 
where $\lfs(\vecy)$ is the  local feature size at point $\vecy\in\Surface$ defined as the minimum Euclidean distance from $\vecy$ to the medial axis $\Medial(\Surface)$~\cite{Amenta1999:Surface}. The medial axis $\Medial(\Surface)$ is defined as the set of points in $\Rd$ where there is more than one closest point. Note that when $\Surface$ is a watertight surface the definition of the medial axis that we use contains both the interior part that is bounded by $\Surface$ and the exterior part that lies outside the bounded domain.

Within the neighborhood $\Tube(\Surface)$, or a subset of it, surface differential operators can be replaced by Cartesian differential operators with \emph{closest point extensions} (see \cite{Ruuth2008:CPM, Marz2012:CPM}). The closest point extension operator $E$ extends surface functions onto $\Tube(\Surface)$ to be constant in the normal direction of $\Surface$ and is defined as $E\sol_\Surface(\vecx) = \sol_\Surface(\cps(\vecx))$. For functions $\sol\in\Tube(\Surface)$ the extension $E$ acts on the restriction of $\sol$ to the surface, i.e., $E\sol = E(\sol|_{\Surface}).$ The Laplace-Beltrami operator in \cref{eq:surface_pde} is equivalent to the following:
\begin{equation}\label{eq:laplace_replacement}
    \mylaplace_\Surface \sol_\Surface(\vecy) = \mylaplace [ E\sol_\Surface](\vecy), \quad \vecy \in \Surface.
\end{equation}
To define the embedding PDE on $\Tube(\Surface)$, we also extend $\source_\Surface$ as $\source(\vecx) = E\source_\Surface(\vecx)$. The equation $\mylaplace [ E\sol_\Surface](\vecx) = \source(\vecx)$, for $\vecx \in \Tube(\Surface)$, is ill-posed because $\source$ is constant in the normal direction of $\Surface$ but $\mylaplace [ E\sol_\Surface]$ is not guaranteed to be. Therefore, the embedding PDE for \cref{eq:surface_pde} becomes
\begin{equation}\label{eq:embedding_pde}
    \mylaplace [ E\sol_\Surface](\vecx) = \source(\vecx) + g(\vecx), \quad \vecx \in \Tube(\Surface),
\end{equation}
where $g(\vecx)$ is a function that compensates for $\mylaplace [ E\sol_\Surface]$ not being constant in the normal direction of $\Surface$. The function $g(\vecx)$ is nonzero for $\vecx \in \Tube(\Surface)\setminus \Surface$ and $g|_\Surface=0$ to ensure \cref{eq:embedding_pde} is consistent with the surface PDE (\cref{eq:surface_pde}) on $\Surface$. Any function $g$ with these conditions has the form $g(\vecx) = \gamma(v(\vecx) - Ev(\vecx))$, where $\gamma \in \mathbb{R}$ and $\gamma \neq 0$.

The Macdonald-Brandman-Ruuth approach (see \cite[Section 2.3]{Chen2015:CPM}) takes $v|_\Surface = u_\Surface$ to allow \cref{eq:embedding_pde} to be written as an equation in one unknown, $v(\vecx)$, since $Eu_\Surface = Ev$ (but importantly $v \neq Ev$ except on $\Surface$). We instead do not restrict $v|_\Surface$ to be $\sol_\Surface$ and interpret $g(\vecx)$ as a modification to the source term $f(\vecx)$, then solve for the unknown solution $\sol(\vecx) = E\sol_\Surface$ in \cref{eq:embedding_pde}. As proven by \citet[Section 3.2]{vonGlehn2013:CPM}, $\sol(\vecx) = E\sol(\vecx)$ since $\sol(\vecx)$ is the extension of a surface function. The property that $\sol(\vecx) = E\sol(x) = \sol(\cps(\vecx))$ allows our projection step during the walk in PWoS, detailed in \cref{sec:method}.

We show through our numerical examples that taking $g(\vecx) = 0$ for all $\vecx \in \Tube(\Surface)$, 
provides qualitatively correct results for graphics applications and quantitatively convergent results in most examples in \cref{sec:results}. However, the choice of $g(\vecx) = 0$ causes \cref{eq:embedding_pde} to be ill-posed as discussed above, and we observe some bias in some convergence studies in \cref{sec:convergence} when $\source$ is complex. Interesting future work would involve constructing a more accurate $g(\vecx)$ function to improve convergence.

In the traditional CPM, one solves the embedding PDE inside a narrow tubular subset of $\Tube(\Surface)$ that is within a constant distance to $\Surface$.
For the typical grid-based variant, the tubular subset is spatially discretized with a grid of uniform spacing. Interpolation and finite differences are applied on the grid to approximate the closest point extension and the spatial Cartesian differential operators, respectively, and then the resulting linear system is solved. Other variants of CPM~\cite{PETRAS2016, PIRET2012, CHEUNG2015194} also require some discretization within $\Tube(\Surface)$.
Thus, while traditional CPM is agnostic to the specific surface representation, it still discretizes the embedding space and solves a globally coupled system.
Moreover, imposing exterior or interior boundary conditions requires tedious grid operations~\cite{King:2023:CPM}.
By contrast, we develop a spatial discretization-free, pointwise Monte Carlo estimator for surface PDEs by incorporating the closest point extension concept into WoS.

When there are Dirichlet boundaries $\Boundary \subset \Surface$, on which the solution $\sol_\Surface$ is given, one can geometrically extend the boundary itself out into the embedding space in the normal directions, assigning it the same boundary value in accordance with the closest point extension. Note that such boundaries need not coincide with the geometric boundaries of the surface itself. 
In the context of grid-based CPM, \citet{King:2023:CPM} discuss how to impose such boundary conditions by duplicating degrees of freedom near the extended boundary.
In our work, we devise a method that uses only the closest point function $\cpc(\vecx)$ to the (pre-extension) boundary $\Boundary$, without the need to construct the extended boundary geometry or perform any complex duplication of degrees of freedom. 

\section{Method}\label{sec:method}
\paragraph{Input} While our algorithm is generalizable to other configurations, we describe our method for the case when $\Surface$ is embedded in $\Rthree$ and $\dim(\Surface) = 1, 2$. Recall that we use the word surface to refer both to ``surfaces'' with $\dim(\Surface)=1$ (curves) as well as $\dim(\Surface)=2$ surfaces. We also allow mixed codimension where parts of the surface have $\dim(\Surface) = 1$ and the rest of the surface has $\dim(\Surface) = 2$. We assume that we can query the closest point function $\cps(\vecx)$ for $\vecx\in\Rthree$. 
Additionally, for surfaces with $\dim(\Surface) = 2$, we assume that we can query the unoriented normal direction of the surface $\vecn(\vecx)$ for $\vecx \in \Surface$. For surfaces with $\dim(\Surface) = 1$, we assume that we can query the tangential direction of the surface $\vect(\vecx)$ for $\vecx \in \Surface$. 
These assumptions are valid for common surface representations, including, but not limited to, polygonal meshes, oriented point clouds, and implicit functions. The theory discussed in \cref{sec:cpm} is based on the assumption that $\Surface$ is smooth; in practice, we observe that applying our technique on discretized surfaces with sharp features behaves well as we show in \cref{sec:convergence}.
Additionally, we assume the Dirichlet boundaries $\Boundary$ have a lower dimension than the dimension of $\Surface$ and support the closest point query $\cpc$. When solving a two-sided boundary value problem for boundaries with $\dim(\Boundary)=1$, we also assume that we can query the tangent direction of $\Boundary$.

\definecolor{myhighlightcolor}{rgb}{0, 0.443, 0.737}
\begin{algorithm}[b]
    \DontPrintSemicolon
    \caption{Projected Walk on Spheres}
    \label{alg:pwos}
    \SetKwInOut{Input}{Input}
        \Input{surface $\Surface$, boundary $\Boundary$, evaluation point $\vecx\in \Surface$, sample walk count $N_P$, volume sample count $N_V$, tolerance $\epsilon$}
    \vskip 0.5EM
    \SetKwFunction{FMain}{EstimateSolution}
    \SetKwProg{Fn}{Function}{:}{}
    \Fn{\FMain{$\Surface$, $\Boundary$, $N_P$, $N_V$, $\vecx$, $\epsilon$}}{
        {\color{myhighlightcolor} $\mathcal{M} \leftarrow$ \texttt{medialAxisPointCloud(}$\Surface$\texttt{)} \tcp*{\cref{sec:lfs}}}
        $\estimate{u}_{\mathrm{sum}} \leftarrow 0$\\
        \For {$n\leftarrow 1$ to $N_P$}{
            $\estimate{u}$ $\leftarrow$ \texttt{RecursiveEstimate(}$\Surface$, $\mathcal{M}$, $\Boundary$, $N_V$, $\vecx$, $\epsilon$\texttt{)} \\
            $\estimate{u}_{\mathrm{sum}}$ $\leftarrow$$\estimate{u}_{\mathrm{sum}} + \estimate{u}$\\
        }
        \Return $\estimate{u}_{\mathrm{sum}}/N_P$
    }
    \vskip 1EM
    \SetKwFunction{FMain}{RecursiveEstimate}
    \Fn{\FMain{$\Surface$,  $\mathcal{M}$, $\Boundary$, $N_V$, $\vecx$, $\epsilon$}}{
        \mbox{\color{myhighlightcolor}$\distance\leftarrow$ \texttt{distanceToBoundary(}$\Surface$, $\mathcal{M}$, $\Boundary$, $\vecx$\texttt{)}\tcp*[r]{\cref{sec:dirichlet}}}\\
        \uIf {$\distance < \epsilon$} {\Return $\sol(\cpc(\vecx))$}
        {\color{myhighlightcolor}$l\leftarrow$ \texttt{localFeatureSize(}$\Surface$, $\mathcal{M}$, $\vecx$\texttt{)}  \tcp*{\cref{sec:lfs}}}
        $\radius\leftarrow \min(l, \distance)$\\
        $\vecy$ $\leftarrow$  \texttt{uniformSphereSample(center=}$\vecx$ , \texttt{radius=}$\radius$\texttt{)}\\
        \mbox{$\estimate\sol_\mathrm{sphere}\leftarrow$ \texttt{RecursiveEstimate(}$\Surface$, $\Medial$, $\Boundary$, $N_V$, {\color{myhighlightcolor}$\cps(\vecy)$}, $\epsilon$ \texttt{)}}\\
        $\{\vecz_1, ... \vecz_{N_V}\}\leftarrow$  \texttt{ballSample(center=}$\vecx$ , \texttt{radius=}$\radius$\texttt{)}\\
        $\estimate\sol_\mathrm{ball}\leftarrow\frac{1}{N_V}\sum_{i=1}^{N_V}\frac{G(\vecx, \vecz_i) \source({\color{myhighlightcolor}\cps(\vecz_i)})}{p(\vecz_i)}$\\
        \Return $\estimate\sol_\mathrm{sphere} + \estimate\sol_\mathrm{ball}$\\
    }
\end{algorithm}

\begin{figure}[b]
\centering
\resizebox{\linewidth}{!}{
\includegraphics{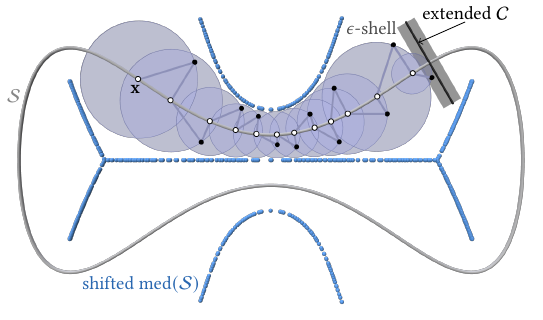}
}
\caption{A PWoS path for the Laplace equation on a gray 1D (curve) surface embedded in 2D space, starting from $\vecx$ and terminating at the extended Dirichlet boundary $\Boundary$.}
\label{fig:pwos}
\end{figure}

\paragraph{Overview} The core idea of our method is to apply the WoS recursive relationship within $\Tube(\Surface)$ while utilizing the closest point extension constraint that $\sol(\vecx) = \sol(\cps(\vecx))$. To do so, we modify the walk process to use spheres contained within $\Tube(\Surface)$ and to \emph{project} the walk position at each recursion step, as illustrated in \cref{fig:pwos}.

The problem we solve is the embedding PDE $\mylaplace u(\vecx) = \source(\vecx) + g(\vecx)$ within $\Tube(\Surface)$. The Monte Carlo estimate of \cref{eq:wos} holds by assuming $g(\vecx)=0$ because the embedding PDE is defined with the Cartesian differential operator.
To estimate the surface PDE's solution at point $\vecx \in \Surface$, we consider a 3D ball centered at $\vecx$ and fully contained inside $\Tube(\Surface)$. Theoretically, it should be the largest ball fully contained inside $\Tube(\Surface)$ that does not cross the extended Dirichlet boundaries $\Boundary$, to minimize the number of steps needed to reach the boundary. We determine the radius of such a ball by taking the minimum of a conservative (under-)estimate of the local feature size at point $\vecx$ (\cref{sec:lfs}) and the distance to the (extended) Dirichlet boundaries (\cref{sec:dirichlet}). 
In \cref{eq:wos}, the Monte Carlo estimate of the solution on the sphere, $\estimate \sol(\vecy)$, needs to be evaluated at point $\vecy$, which does not lie on $\Surface$ in general.
The closest point extension constraint of $\sol$ provides a convenient relationship here: 
the embedding PDE's solution at point $\vecy$ should coincide with the surface PDE's solution at the projected point, $\cps(\vecy)$.
We can therefore project the point $\vecy$ onto $\Surface$ at the end of each recursion step hence $\estimate u(\vecy) = \estimate u(\cps(\vecy))$.
After this projection at the end of each step, we continue the recursion.
The source term similarly uses the closest point projection for the closest point extension, replacing the recursive relationship of WoS (\cref{eq:wos}) with\looseness=-1
\begin{equation}\label{eq:pwos}
\estimate\sol(\vecx) = \estimate\sol(\cps(\vecy)) + \frac{1}{N_V}\sum_{i=1}^{N_V}\frac{G(\vecx, \vecz_i) \source(\cps(\vecz_i))}{p(\vecz_i)}.
\end{equation}
We choose $p(\vecz_i)\propto 1/\lVert \vecx- \vecz_i \rVert_2$ in our implementation.

Analogous to the original WoS method, we terminate the recursion when the point $\vecx$ falls within a distance $\epsilon$ of the (extended) Dirichlet boundary by taking the boundary value $\sol(\cps(\vecx))$. We provide pseudocode for an instance of our algorithm in \cref{alg:pwos}, where we highlight the difference between our proposed method and WoS. We also provide a visualization of a potential path of our algorithm when $\Surface$ is embedded in $\Rtwo$ in \cref{fig:pwos}. Notably, PWoS is a generalization of the WoS algorithm: when $\dim(\Surface) = \dimension$ (i.e., the codimension-zero case), the local feature size is infinite, the distance to the Dirichlet boundary $\Boundary$ can be computed with the closest point query $\cpc$, and the last closest point projection of $\vecy$ has no effect since $\cps(\vecy) = \vecy$.
When $\dim(\Surface) < \dimension$, in addition to the closest point projection at the end of each recursion step, our algorithm utilizes two nontrivial steps: the local feature size estimation using a medial axis point cloud and the computation of the distance to the (extended) Dirichlet boundary. We discuss these in the following two subsections.\looseness=-1

\begin{figure}[b]
\centering
\begin{minipage}[hb]{0.3\linewidth}
\begin{center}
\includegraphics[trim={9cm 1cm 9cm 0.5cm}, clip, width=\linewidth]{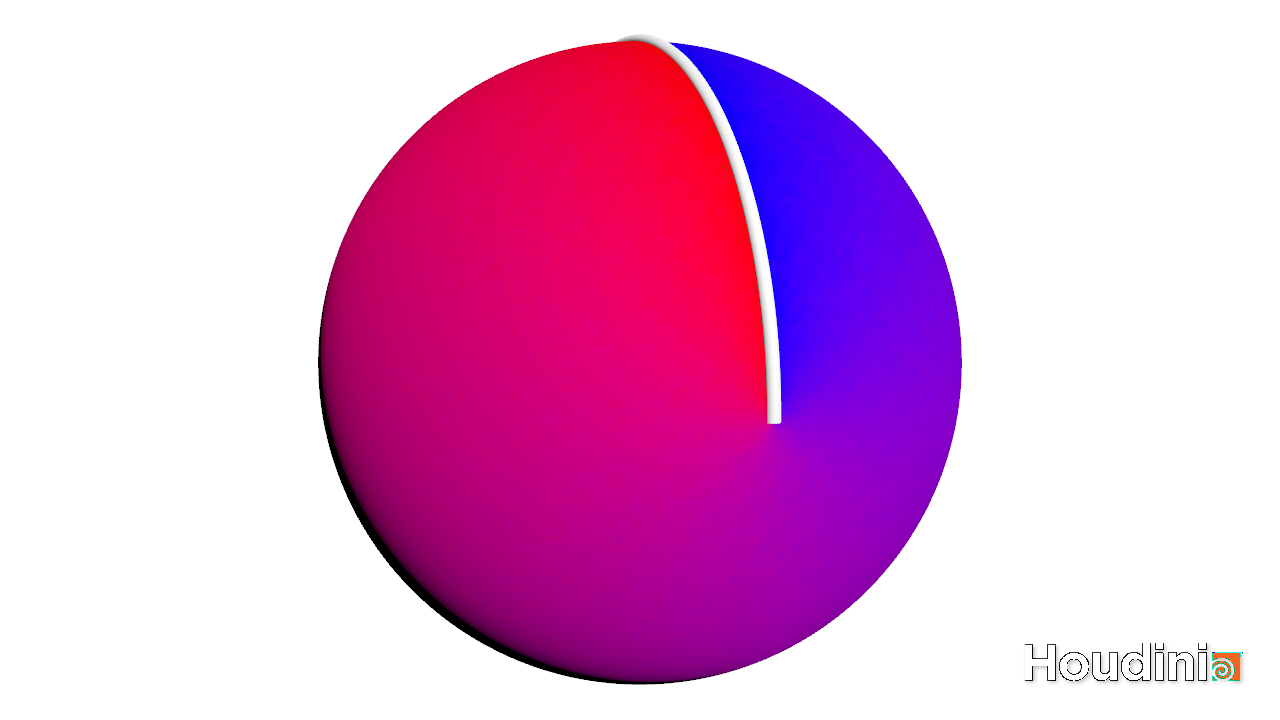}
\end{center}
\end{minipage}
\begin{minipage}[hb]{0.6\linewidth}
\begin{center}
\begin{tabular}{|c|c|c|}\hline
        Local Feature Size        &  Avg. Step Count  \\\hline
        0.99 &31.1398\\
        0.5 &47.84\\
        0.25 &128.981\\
        0.125& 462.781\\
        0.0625 &1818.73 \\\hline          
\end{tabular}
\end{center}
\end{minipage}
\caption{Average number of steps required with different conservative local feature size estimates. While any positive value smaller than 1 is valid for this setup, using a local feature size estimate that is too small leads to excessively long walks.}
\label{fig:badlfsestimate}
\end{figure}

\subsection{Local Feature Size Estimation}\label{sec:lfs}
To determine the radius of the sphere centered at $\vecx\in\Surface$ that is fully contained inside $\Tube(\Surface)$ at each recursion step of the walk, we need a conservative lower bound estimate for the local feature size at $\vecx$.
One naive approach would be to use a small enough positive constant value for all $\vecx\in\Surface$, similar to the grid-based CPM ~\cite{Ruuth2008:CPM}. This is a valid strategy, but it would often yield a sphere radius smaller than necessary, requiring more recursion steps for the walks to reach a Dirichlet boundary and making the method inefficient. \cref{fig:badlfsestimate} illustrates a result of our algorithm on a unit sphere using different (artificial) local feature size estimates. The analytical local feature size of this surface is $1$ everywhere. Although using any smaller value would still give consistent results, it significantly increases the average step count for the walks. For more complex shapes, one cannot compute the analytical local feature size in general and using a small constant value in its place is inefficient. This issue motivates our need for a better local feature size estimate.\looseness=-1

To estimate the local feature size, we compute a point cloud approximation of the medial axis and estimate the local feature size as the distance from any query point $\vecx \in \Surface$ to its nearest point in the medial axis point cloud.
One could use any local feature size estimation algorithm and/or medial axis extraction algorithm (see e.g., \cite{Tagliasacchi2016}), such as one that outputs or uses the medial axis' connectivity.
Since such methods are often comparatively costly, we employed a simple point cloud-based method, which we describe below.

\paragraph{Medial ball extraction}
We first densely scatter points $\vecx_i$ inside a ball in $\Rthree$ having a radius equal to half the length of the diagonal of the bounding box of $\Surface$, so the entire surface is fully enclosed.
For each point $\vecx_i$, we use the closest point query $\cps(\vecx_i)$ to compute its two opposing normal directions at $\cps(\vecx_i)$.
Specifically, we normalize the vector $\vecx_i-\cps(\vecx_i)$ to get the first direction and invert its direction to get the second one. We then use the method of \citet{Ma2012Medial} to extract a point cloud that represents the medial axis, as follows. For each side of the surface (i.e., each normal), we start with a large sphere tangent to $\cps(\vecx_i)$, whose center necessarily lies on the normal ray. The initial radius of the ball is set to the length of the diagonal of the bounding box of $\Surface$. Then, we iteratively shrink the size of the sphere, moving its center to maintain tangency at the surface point $\cps(\vecx_i)$, until the closest surface point from the center of the sphere does not change. This algorithm gives \emph{medial balls}, i.e., balls with their centers on the medial axis.
As the number of initial scattered points increases, the extracted point cloud balls tend to cover the entire medial axis.
While \citet{Ma2012Medial} assumed that the surface is represented as an oriented point cloud, we observed that the algorithm works well with other surface representations by adjusting its termination criteria. See the paper by \citet{Ma2012Medial} and our implementation for details.

\paragraph{Scale axis pruning}
Directly using the medial ball centers as the medial axis point cloud does not work well in general when $\Surface$ contains any noise or artificial sharp corners introduced by the discretization of a smooth surface. We therefore prefer to estimate (only) the \emph{stable} part of the medial axis, which is not affected by any small perturbation of $\Surface$. A common solution is, therefore, to \emph{prune} unstable components of the medial axis, which itself remains an active research topic~\cite{Tagliasacchi2016}. We take inspiration from the scale axis transform (SAT)~\cite{Giesen2009, Miklos2010SAT}, but design an alternative since SAT considers topology information of the medial axis, which is unnecessary for local feature size estimation. Our alternative is simpler and faster since topology information is omitted.
We first scale all the medial balls by a constant factor $s>1$. Then, for any pair of medial balls $B_{r_1}(\vecx_1)$ and $B_{r_2}(\vecx_2)$, we remove the smaller ball from the set of medial balls if it is fully contained in the other. That is, if $s \cdot r_1 < s \cdot r_2 + \lVert \vecx_1 - \vecx_2\rVert_2$, we remove the ball $B_{r_1}(\vecx_1)$.

Note that some of the medial balls may have a very large radius before pruning. For example, an exterior medial ball for a surface of a convex shape would have an infinite radius in theory (but our algorithm returns at most the length of the diagonal of the bounding box of $\Surface$).
When such large medial balls are used in our pruning algorithm, they can easily (and undesirably) remove important, stable parts of the medial axis. The original SAT approach was applied only to the interior medial axis of closed surfaces. Therefore, this issue was not observed since the interior medial ball radii are always bounded and proportional to the size of local features of $\Surface$.
To address this problem, we consider each pair of tangent balls generated at the same surface point, and replace the larger one with a tangent ball having the radius of the smaller one. In other words, before we prune the medial axis, we shift the medial ball centers of the larger medial ball in the pair and shrink its radius. In \cref{fig:pwos}, we visualize the medial ball centers after this shifting operation.

After this pruning, the set of medial ball centers gives the medial axis point cloud we use to estimate the local feature size.
Adjusting the scaling parameter $s$ allows us to control the pruning strength. Unless otherwise stated, we use the value $s=1.15$ for our results.

\paragraph{Conservative and nonzero local feature size}
As the medial axis is represented as a point cloud and the nearest point distance may give a larger value than the actual local feature size, we multiply by a small constant ($0.9$ in our implementation) to ensure a conservative estimate of the local feature size. 
When there are sharp corners in the geometry, the analytical local feature size is zero, and the walk will become stuck. To prevent this problem, when the estimated local feature size is smaller than a positive constant threshold $\lambda$, we return $\lambda$ as the local feature size estimate instead. This process essentially rounds sharp corners with rounding radius $\lambda$. The uniform grid size adopted in grid-based CPM has a similar effect. We do not observe any significant error due to this rounding, as we show in \cref{sec:convergence}.

\subsection{Distance to Extended Dirichlet Boundaries}\label{sec:dirichlet}
Dirichlet boundaries $\Boundary$ are extended in the normal direction of $\Surface$ and the solution in the embedding space on this extended boundary is determined by the closest point extension of Dirichlet values on $\Boundary$.
Therefore, we need to compute the minimum distance to the extended Dirichlet boundaries, and limit the sphere radius in PWoS further if it is less than the local feature size. 
To determine the distance to the extended Dirichlet boundary from point $\vecx \in \Surface$, we first find the closest point that lies on the boundary before the extension, $\cpc(\vecx)$.
The subset of the extended boundary that is extended from $\cpc(\vecx)$ takes the shape of a line segment when $\dim(\Surface)=2$ and a disk when $\dim(\Surface)=1$.
We set the line segment's half-length or the disk radius to the local feature size at $\cpc(\vecx)$ using the algorithm in \cref{sec:lfs}.
We can compute the distance from the point $\vecx \in \Surface$ to this line segment or disk without explicitly constructing the extended boundary geometry. 
When $\dim(\Surface)=2$, the distance $\distance$ to the extended boundary is given by
\begin{equation}
\begin{split}
    \vecr &= \cpc(\vecx) - \vecx,\\
    \distance &= \lVert \vecr - \mathrm{clamp}(\vecr\cdot\vecn, -l,l) \cdot \vecn \rVert_2,
\end{split}
\end{equation}
where $\vecn$ and $l$ are the surface normal and the local feature size estimate at $\cpc(\vecx)$, respectively.
When $\dim(\Surface)=1$, the normal direction is not uniquely defined, so we instead use a similar method based on the surface tangent $\vect$:
\begin{equation}
\begin{split}
    \vecq &= \vecr - (\vecr\cdot\vect)\vect,\\
    \distance &= \begin{cases}
       \lvert \vecr\cdot \vect \rvert, & \mathrm{if}\; \lVert \vecq \rVert_2 < l, \\
       \lVert \vecr -  l \cdot(\vecq /\lVert \vecq \rVert_2) \rVert_2, &\mathrm{otherwise}.
    \end{cases}
\end{split}
\end{equation}

\subsection{Generalizations}
\subsubsection{Screened Poisson Equation} We have so far considered only the Poisson equation. For WoS, \citet{Sawhney:2020:MCGP} proposed a generalization of WoS to the screened Poisson equation $\mylaplace \sol - \sigma \sol = \source$, where $\sigma$ is a positive constant. The embedding PDE for the screened Poisson equation is constructed similarly to \cref{eq:embedding_pde} using closest point extensions~\cite[Section 2.3]{Chen2015:CPM}. Therefore, similar to WoS~\cite{Sawhney:2020:MCGP}, we replace the integral equation (\cref{eq:integraleq}) used in our recursive relationship with
\begin{equation}\label{eq:integraleq_sp}
    \sol(\vecx) = \frac{c_{r, \sigma}}{\lvert \partial B_\radius(\vecx) \rvert} \int_{\partial B_\radius(\vecx)} \sol(\vecy) \dAy + \int_{B_\radius(\vecx)}  \source(\vecz) G_\sigma(\vecx, \vecz) \dVz,
\end{equation}
where $G_\sigma$ is the Yukawa potential and $c_{r, \sigma}$ is a positive number smaller than $1$. To evaluate the first term, instead of multiplying the contribution from the recursion by $c_{r, \sigma}$ as suggested by \citet{Sawhney:2020:MCGP}, we use a Russian Roulette strategy: we terminate the path early with probability $1 - c_{r, \sigma}$ with zero contribution and otherwise we use the estimate of the solution without multiplying  by $c_{r, \sigma}$. As PWoS is a generalization of WoS, we can apply this Russian Roulette strategy to WoS as well.
This scheme allows us to terminate the PWoS path early without waiting for it to reach the boundary and without introducing additional bias. It also makes it possible to apply PWoS to problems without Dirichlet boundaries. We can even use this estimator for the screened Poisson equation as a Tikhonov regularization of the Poisson equation without boundaries: solving the screened Poisson equation with small $\sigma$ yields an approximate solution to the Poisson equation. Similar regularization ideas appear in multiple contexts~[\citetalias{SabelfeldSimonov1994}~\citeyear{SabelfeldSimonov1994}, Section 6.3; \citetalias{sawhney2023walk}~\citeyear{sawhney2023walk}].

\subsubsection{Divergence Source Term and Gradient Estimation}\label{sec:graddiv}
Many graphics applications, such as the heat method for geodesic distance computation~\cite{Crane2013} and the projection step of fluid simulation~\cite{FOSTER1996:Realistic, Stam1999}, give rise to a Poisson equation with a source term expressed as the divergence of a vector field, $\source_\Surface = \grad_\Surface\cdot\vech_\Surface$, where $\vech_\Surface$ is a vector field defined over the surface.
These applications also require the estimation of the gradient of the solution instead of the solution itself. With grid-based CPM, the differential operators defined in the embedding Cartesian domain can be used to solve such problems~\cite{Auer2012, King:2023:CPM}. For our PWoS, however, we do not use any embedding grid structure, and we do not assume any specific surface representation, so we cannot use such discrete differential operators.

To solve a problem with a divergence of vector field as the source term, \citet{Sugimoto:2024:VelMCFluids} proposed to use integration by parts to convert the volume integral arising from the source term to a form that does not explicitly depend on the divergence of the vector field. This was done in the context of the walk on boundary method~\cite{Sugimoto2023:WoB}, which is a Monte Carlo volumetric PDE solver based on a different integral equation formulation, and we adapt this approach to (projected) WoS.
To estimate the gradient with PWoS, we can use the gradient estimator for WoS~\cite[Section 3.1]{Sawhney:2020:MCGP}, because the solution's gradient is zero in the surface normal directions due to the closest point extension constraint $\sol(\vecx) = \sol(\cps(\vecx))$. The gradient estimator replaces the first step of recursion based on \cref{eq:integraleq} with  one we can derive by taking the gradient of \cref{eq:integraleq}. 
We discuss the details of these estimators in the supplemental note and demonstrate an application for the geodesic distance computation in \cref{sec:geodesic}.

\section{Mean Value Filtering with Discrete Basis}
The method we described in \cref{sec:method} is suitable for evaluating the solution at a few evaluation points independently. To improve the efficiency of the method for many evaluation points (i.e., mesh vertices), it is often desirable to utilize the spatial consistency of the solution. For WoS, a few such methods 
have been proposed, but those approaches are not trivially applicable to our setup, where we have additional closest point extension constraints between the surface and embedding PDEs.
Specifically, the boundary value caching approach~\mycite{Miller2023} is based on a boundary integral equation defined in the Cartesian coordinates and is not applicable to surface PDEs. Adaptation of reverse WoS~\cite{Qi2022} or mean value caching~\cite{Bakbouk2023:MVC} to our setting would require a mapping of a PDF on the surface to a PDF on the spheres used in our walk process, and we found it difficult to design such an algorithm for general surfaces.

Inspired by the filtering method of \citet{Bakbouk2023:MVC}, we designed a filtering method that uses discrete basis functions defined over a surface represented as a polygonal mesh or point cloud. When we apply \cref{eq:pwos} at an evaluation point, if $\estimate \sol(\cps(\vecy))$ is a precomputed Monte Carlo estimate, we do not need to continue the recursion and can simply use the precomputed estimate in its place. In general, we do not have the estimate $\estimate \sol(\cps(\vecy))$ available at $\cps(\vecy)$. %
Thus, we get the estimate by interpolating the estimate of solutions already computed at a discrete set of points by accepting the bias due to interpolation. For example, for a triangle mesh, we use barycentric coordinates as the interpolation basis.
This can be interpreted as a PDE-aware smoothing filter of the solution when we consider this process as a weighted average of the solution estimates at nearby evaluation points.

To estimate the solution at all mesh vertices or all points in a point cloud, we first compute an approximate solution to the problem with a low sample count using the method described in \cref{sec:method} and apply this filtering step to get an updated estimate. We can also precompute the filtering weights per evaluation point first and apply the same filter a few times to achieve an even smoother solution without too much additional cost. We can similarly design a gradient estimation filter based on \cref{sec:graddiv}. While this filtering approach utilizes a discrete basis defined over the surface, we do not use any explicit discrete differential operators and do not need to solve a linear system, which is in contrast to the traditional methods based on discrete Laplacians.

The method of \citet{Bakbouk2023:MVC} similarly designed a smoothing filter for volumetric PDEs. %
Their method can keep the estimate unbiased by assuming that the evaluation points are sampled with a known PDF and that the original estimates are unbiased, but our method introduces bias due to the use of a discrete basis. We nevertheless observe that, within a reasonable runtime budget,
our biased filtering method can reduce the error compared to the PWoS algorithm without filtering. We leave the development of an unbiased variant to future work.\looseness=-1

\section{Results}\label{sec:results}
We implemented PWoS in Houdini 20.0~\cite{Houdini} without GPU acceleration, using its built-in closest point queries. Our implementation is provided as supplemental material.

\definecolor{mycolor1}{HTML}{6A99D0}
\definecolor{mycolor2}{HTML}{DE8344}
\definecolor{mycolor3}{HTML}{7EAB55}
\definecolor{mycolor4}{HTML}{F5C342}

\newcommand{\tableelemwidth}{0.115\linewidth}
\newcommand{\addgraph}[3]{
\begin{minipage}[t]{\tableelemwidth}
\includegraphics[trim={10.5cm 0cm 10.5cm 0cm}, clip, width=0.8\linewidth]{#1.png}\\
\hspace*{-2em}
\resizebox{1.2\linewidth}{!}{\includegraphics{main-figure#3.pdf}}
\centering
\textsf{\small #2}
\end{minipage}
}

\begin{figure*}[h]
    \centering

    \addgraph{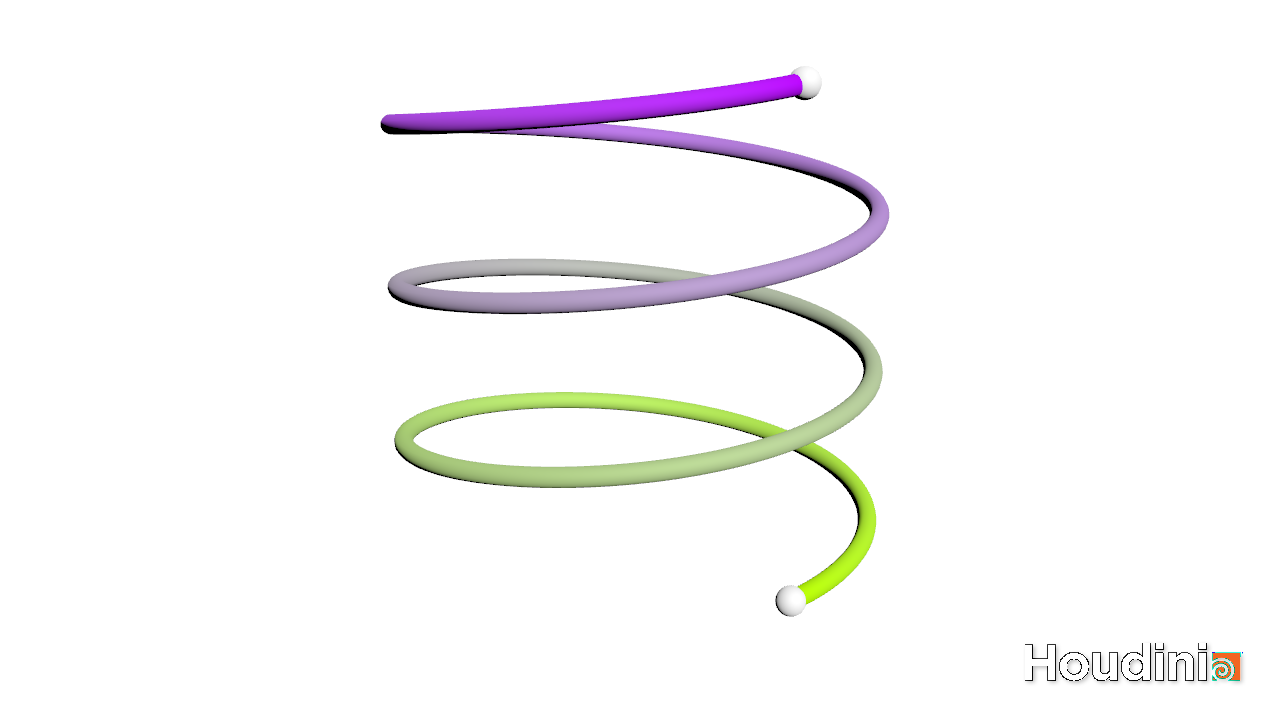}{(a) Helix\\Ends\\Laplace}{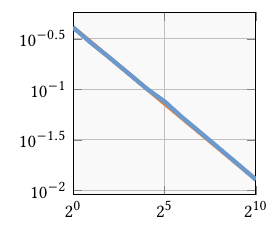}%
    \addgraph{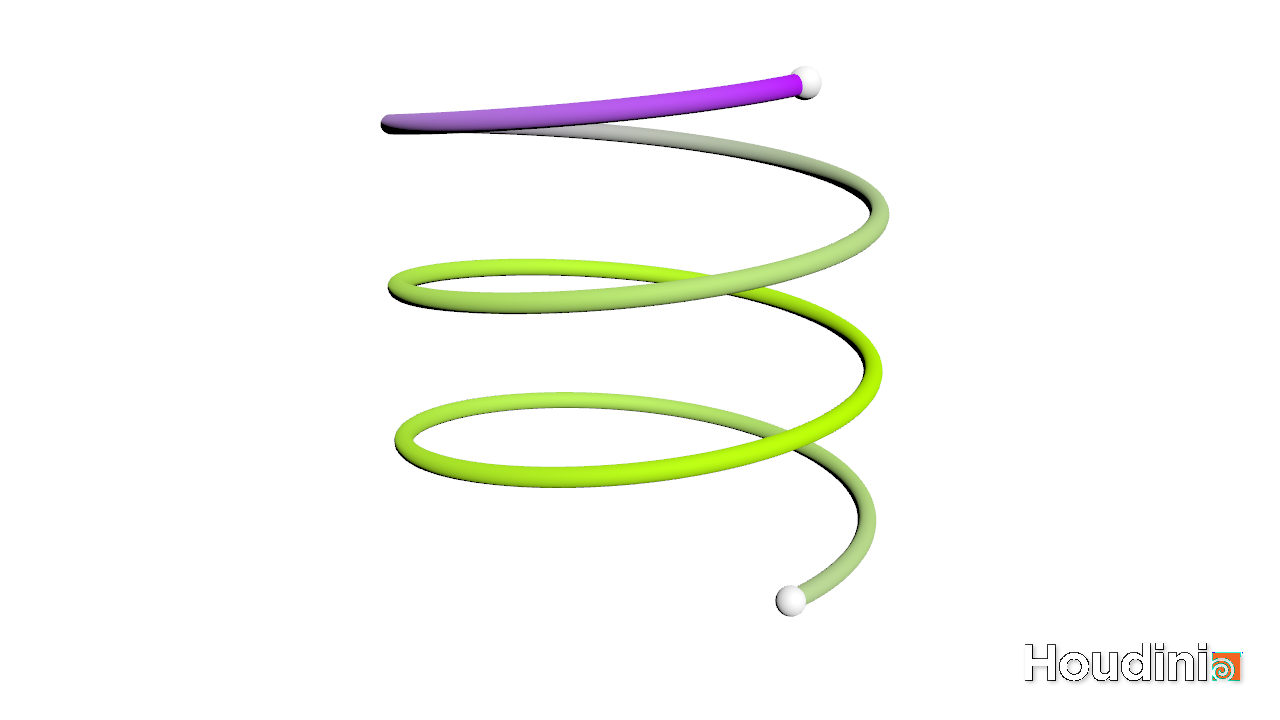}{(b) Helix\\Ends\\Poisson}{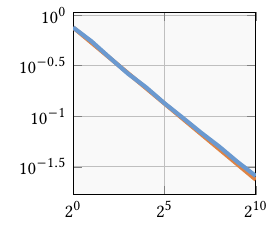}%
    \addgraph{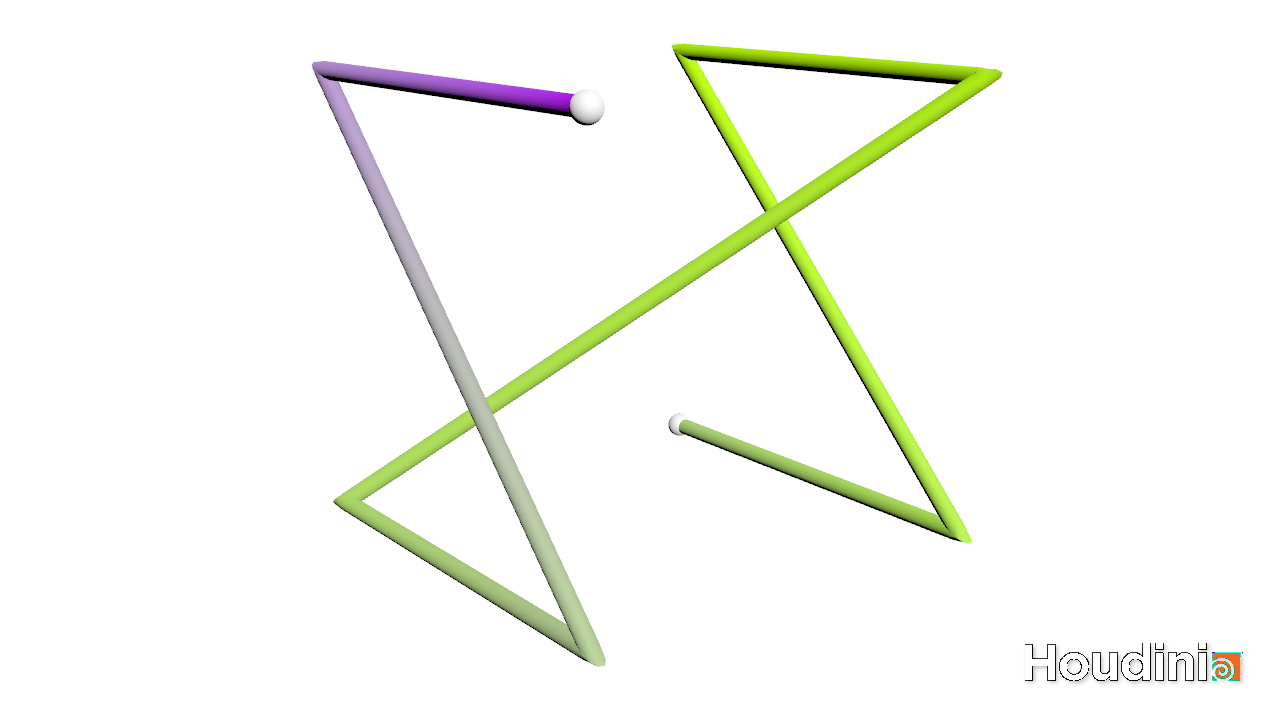}{(c) Z-Order ($s=1.15$)\\Ends\\Poisson}{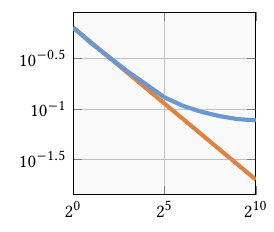}%
    \addgraph{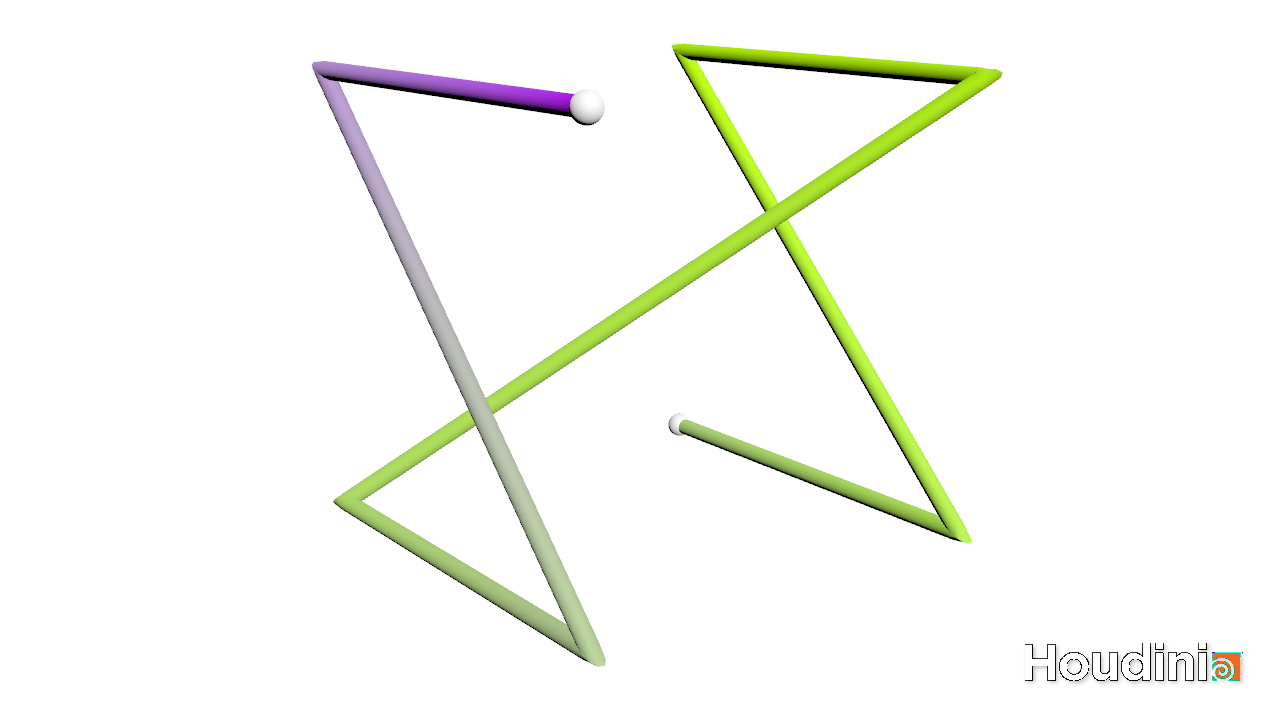}{(d) Z-Order ($s=1.05$)\\Ends\\Poisson}{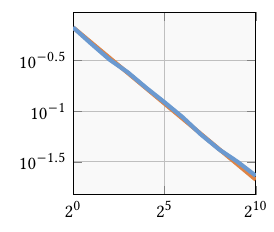}%
    \addgraph{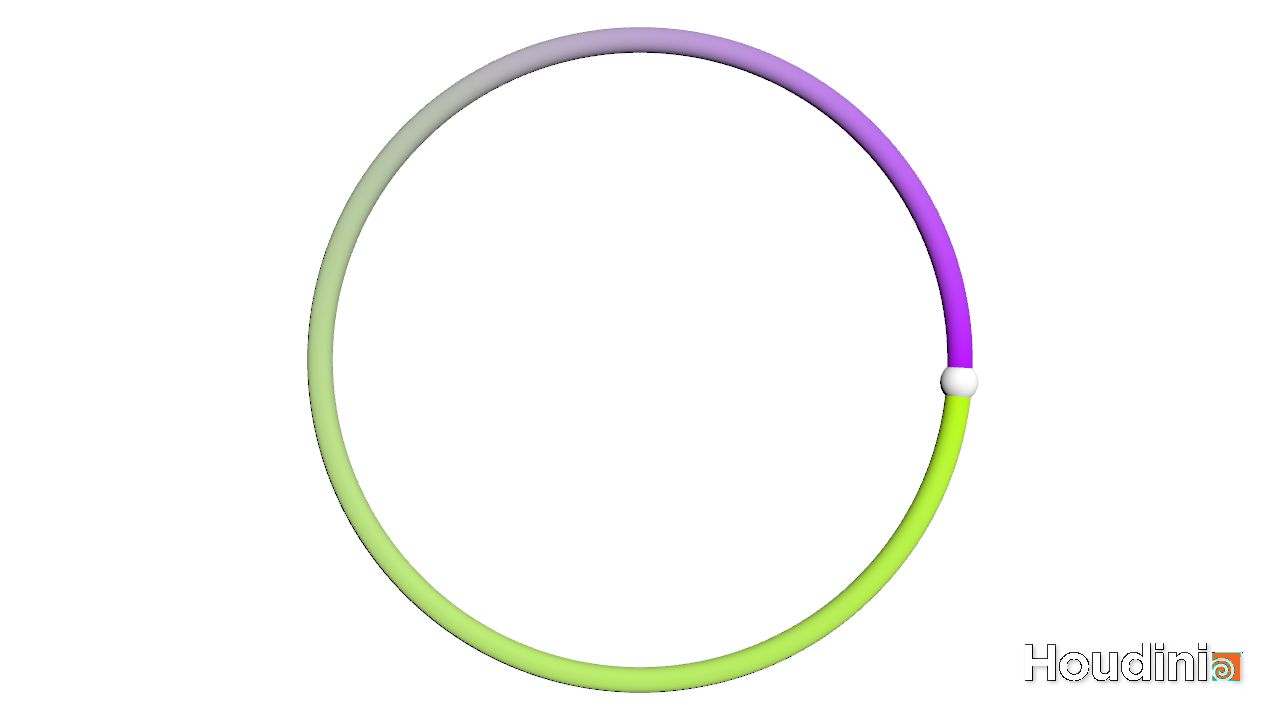}{(e) Circle\\Two-sided\\Poisson}{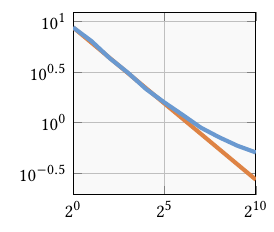}%
    \addgraph{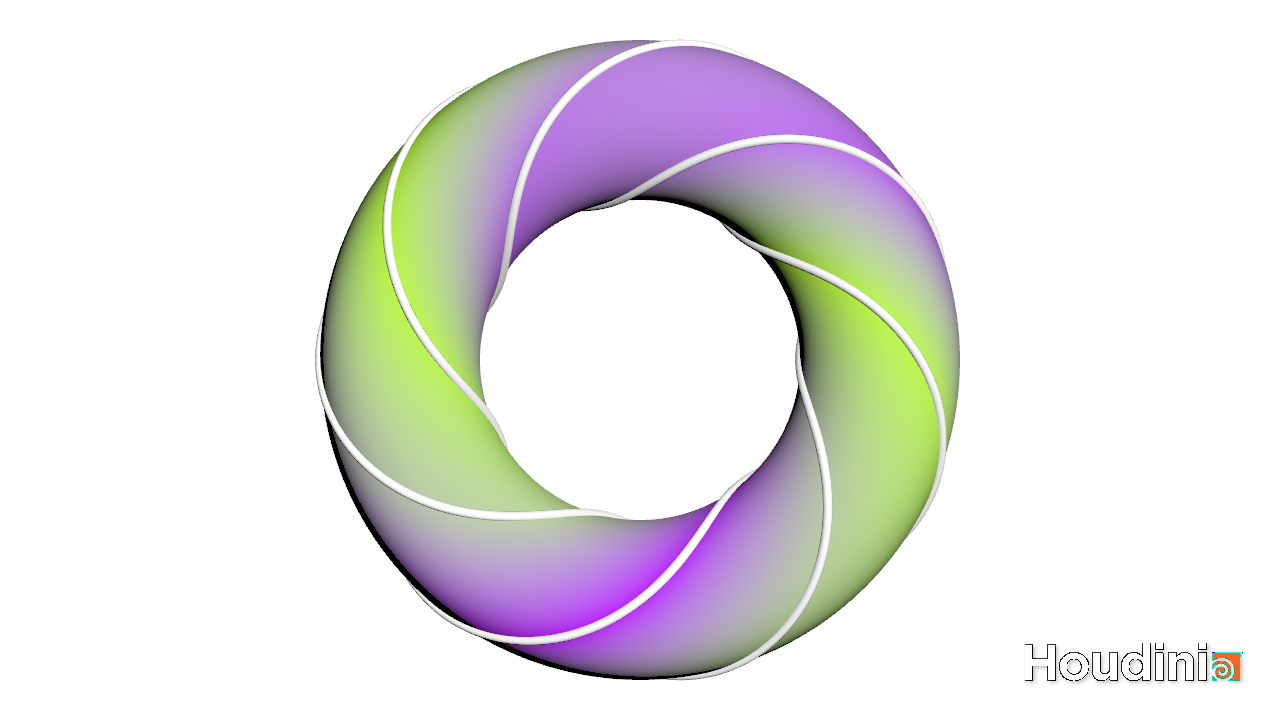}{(f) Torus\\Torus Knot\\Laplace}{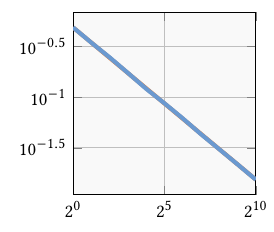}%
    \addgraph{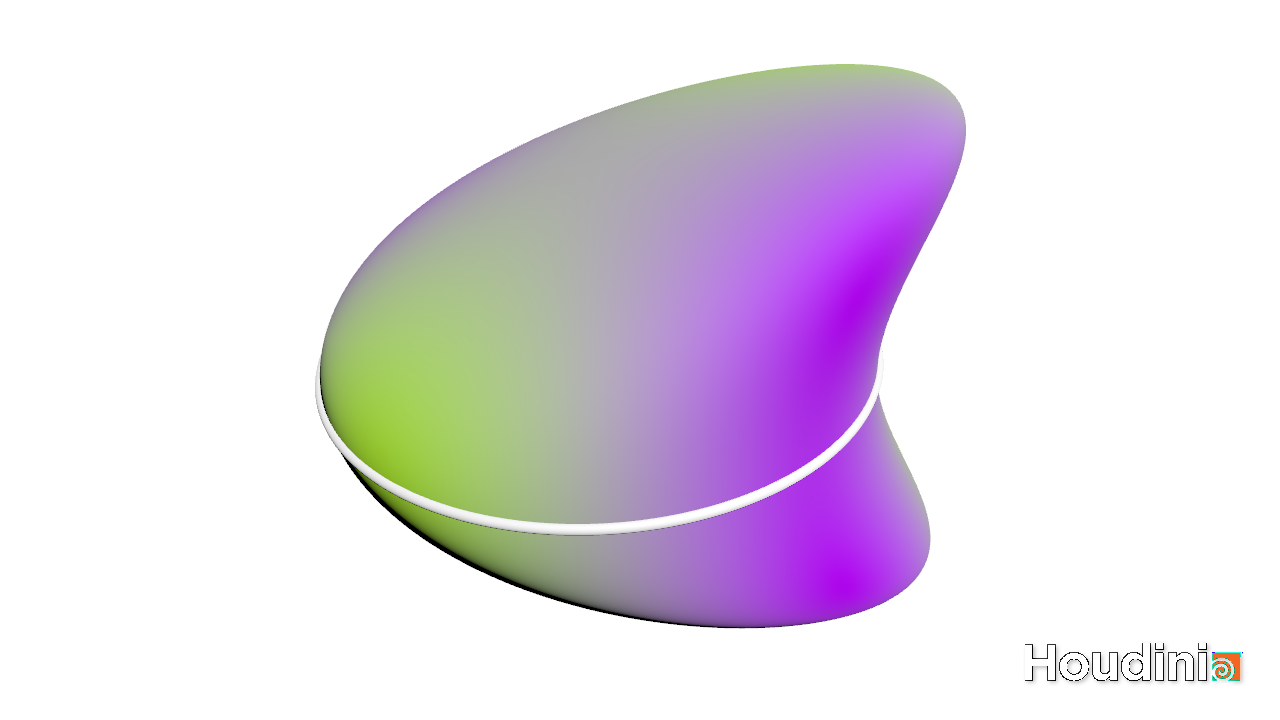}{(g) Dziuk Surface\\Closed\\Screened Poisson}{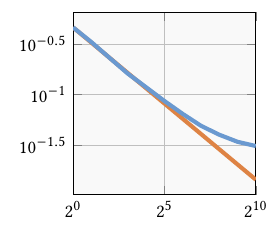}%
    \addgraph{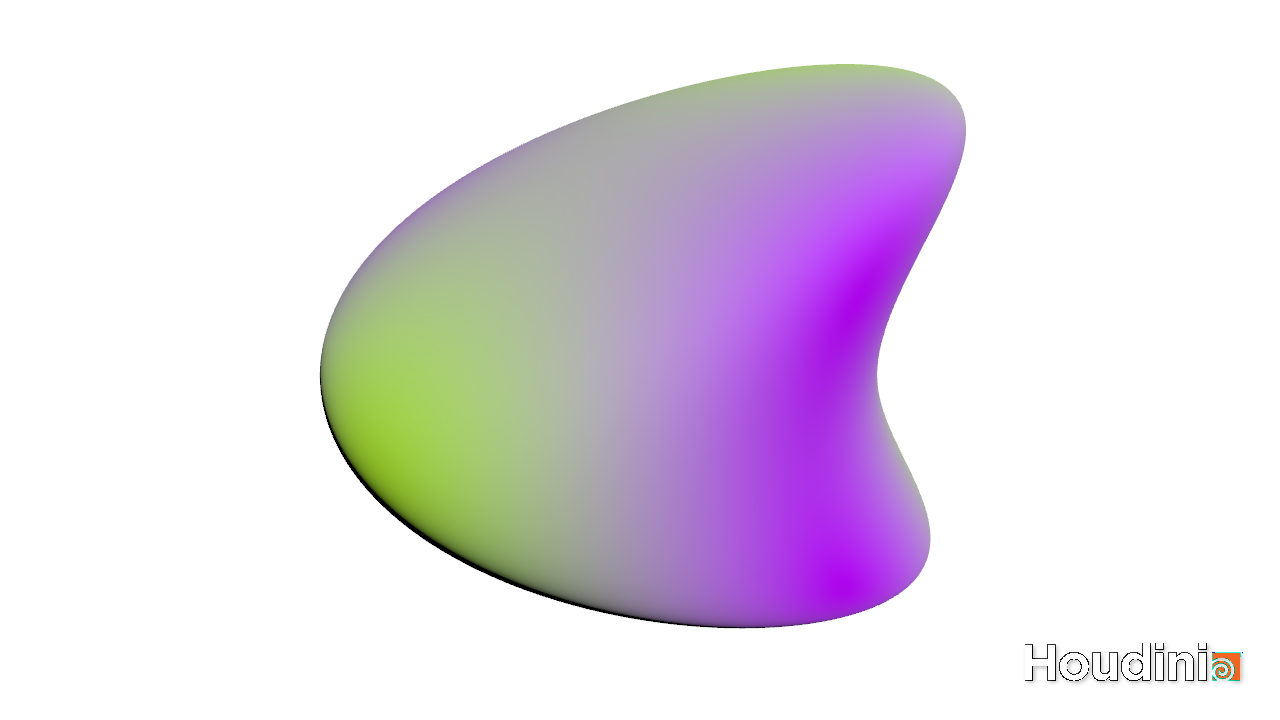}{(h) Dziuk Surface\\No Boundary\\Screened Poisson}{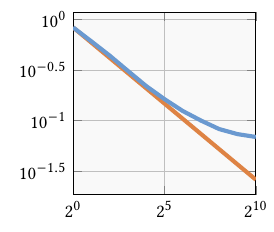}
    
    \addgraph{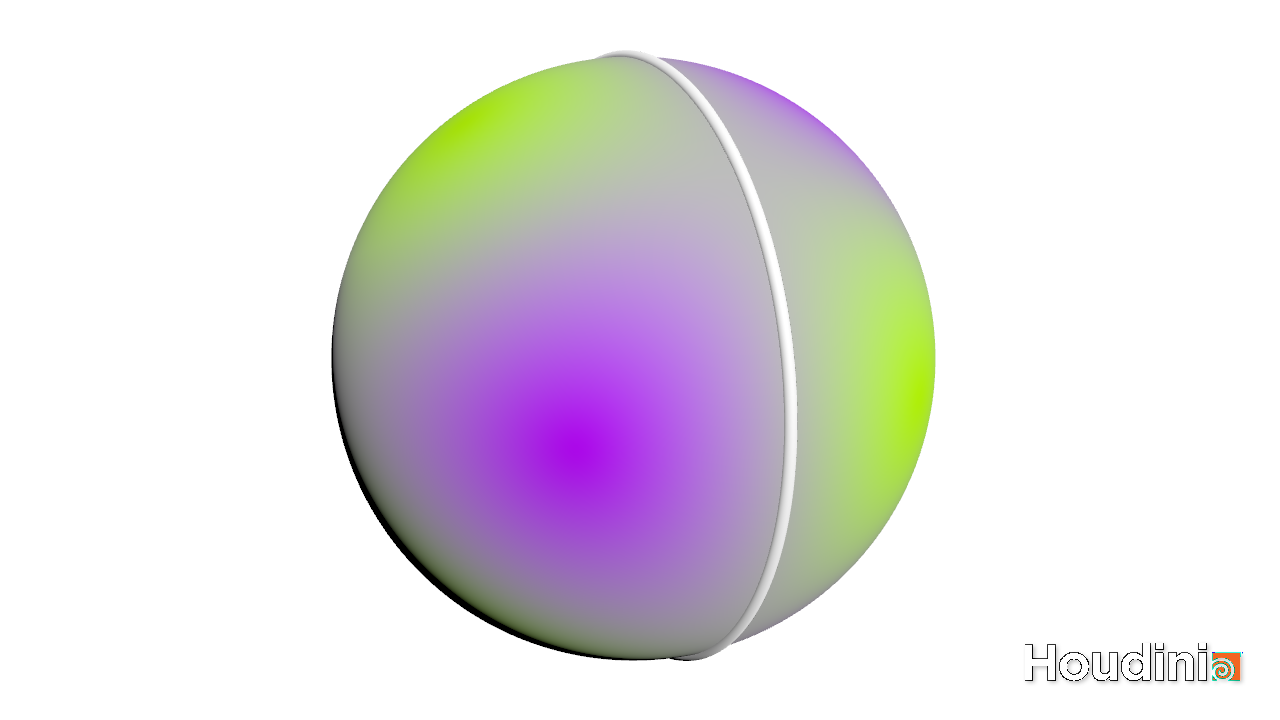}{(i) Sphere\\Closed\\Poisson}{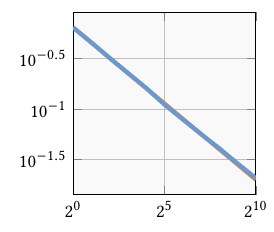}%
    \addgraph{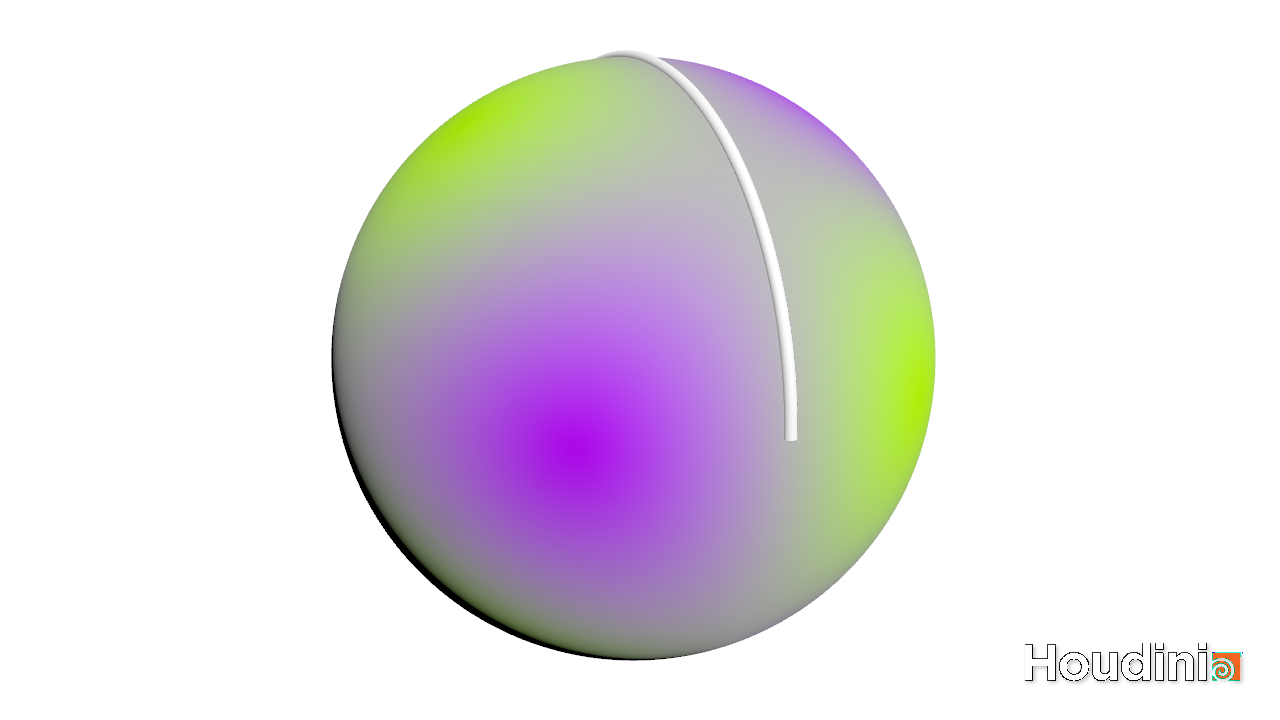}{(j) Sphere\\Open\\Poisson}{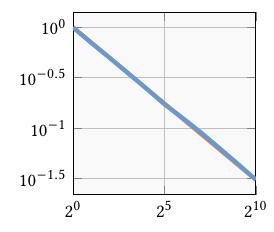}%
    \addgraph{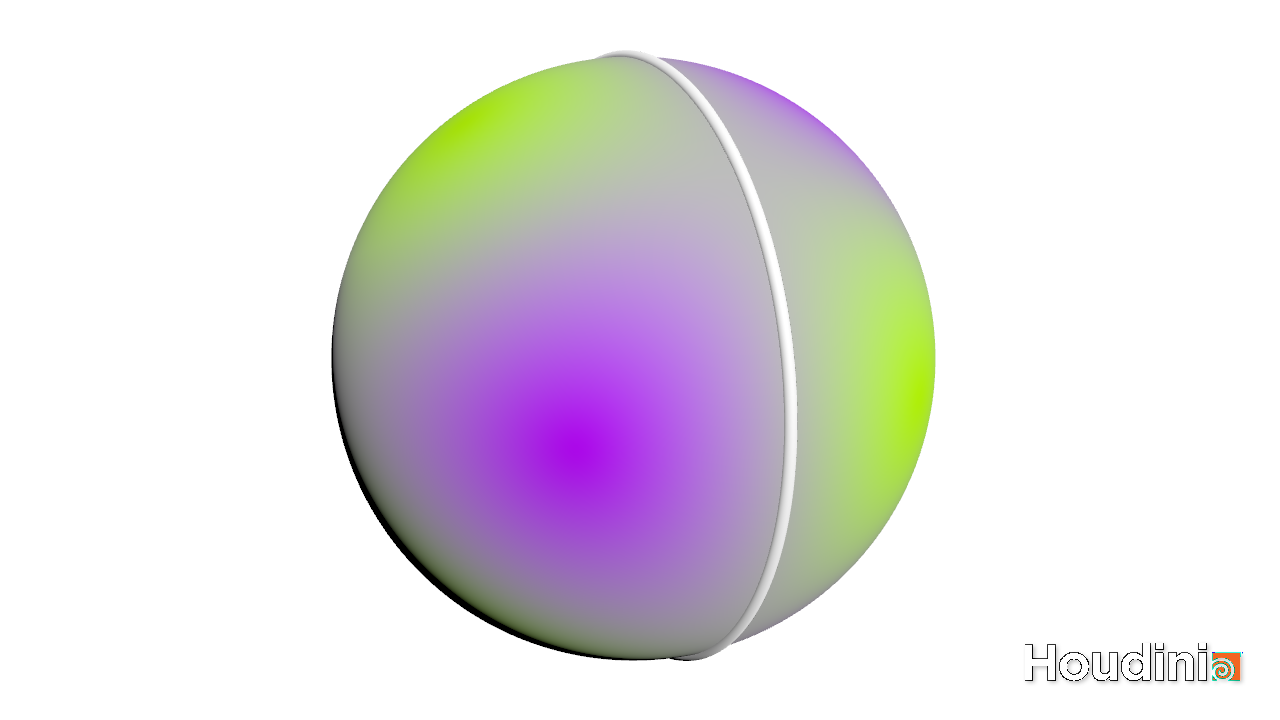}{(k) Sphere\\Closed\\Screened Poisson}{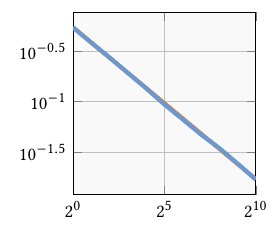}%
    \addgraph{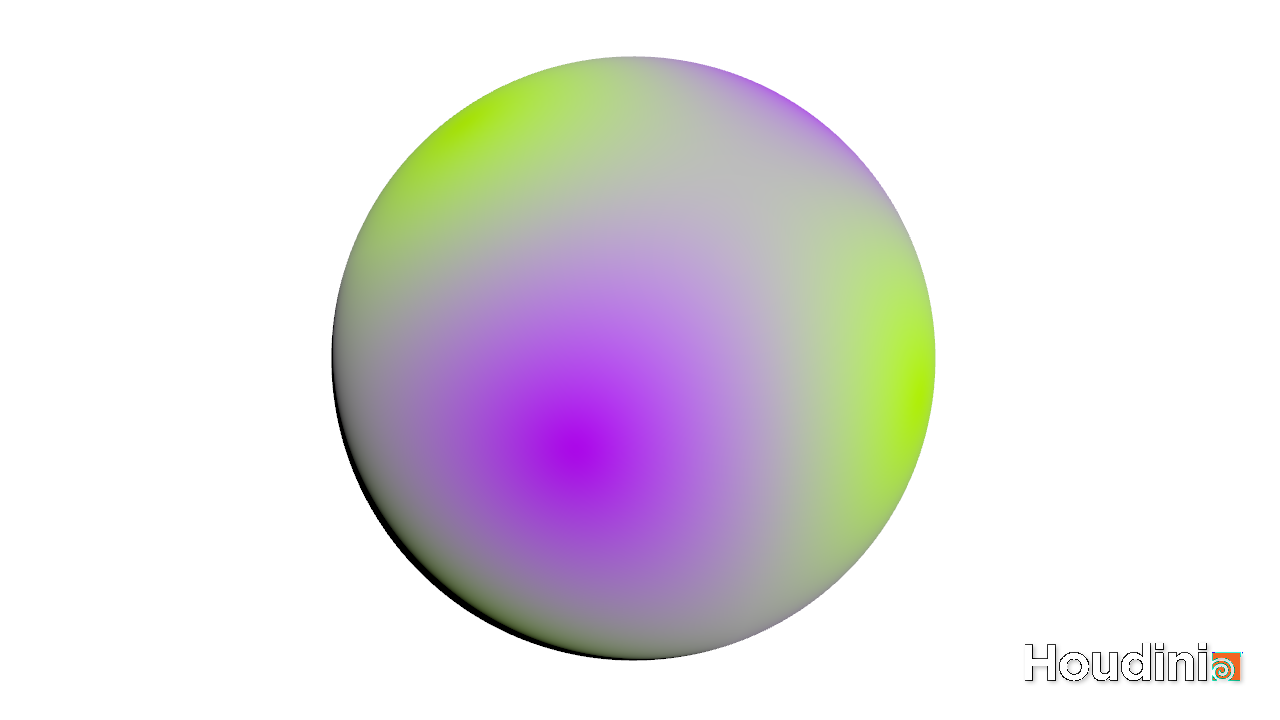}{(l) Sphere\\No Boundary\\Screened Poisson}{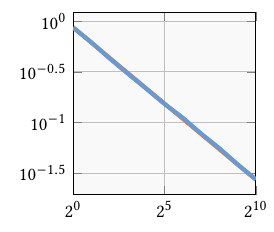}%
    \addgraph{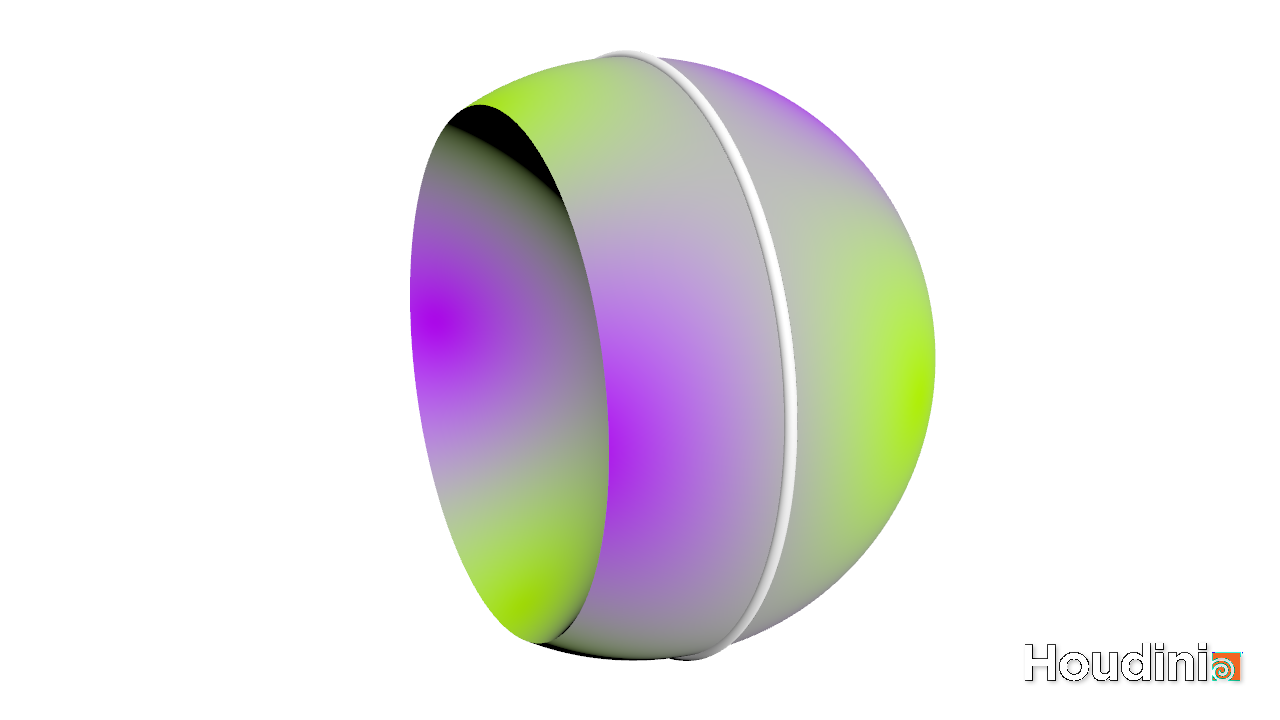}{(m) Punched\\Closed\\Poisson}{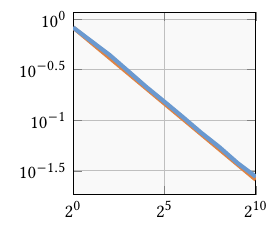}%
    \addgraph{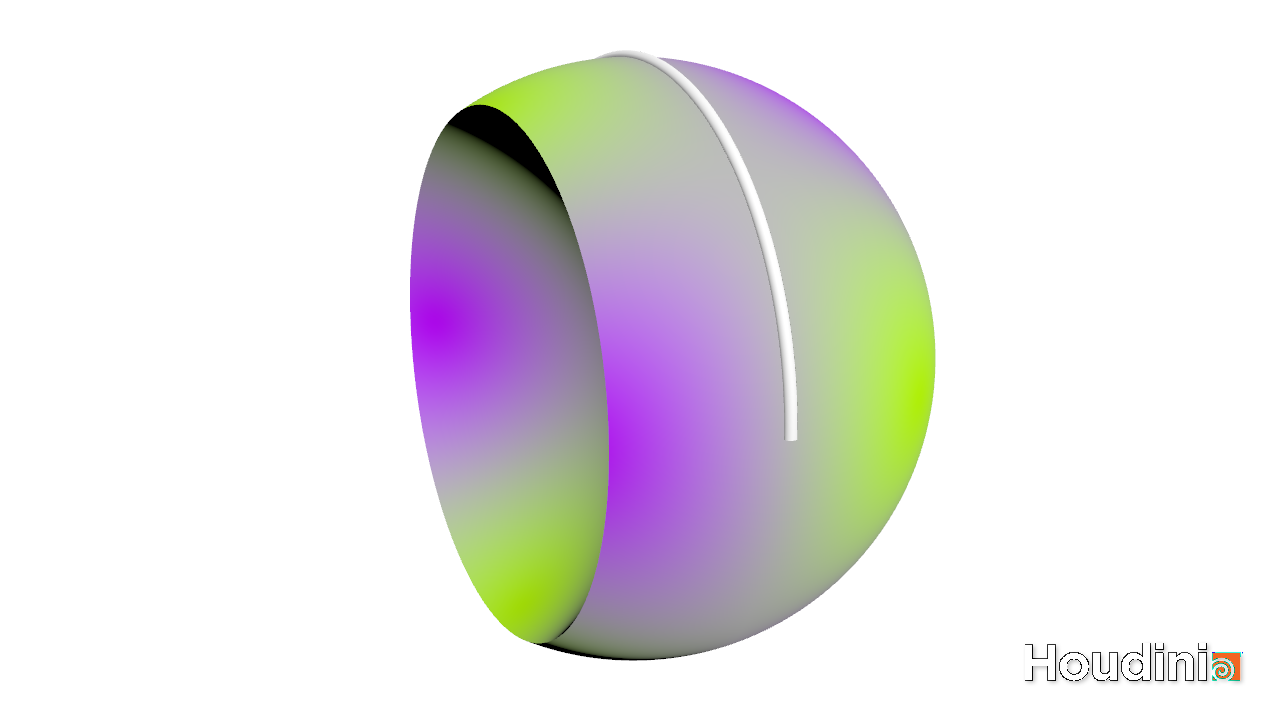}{(n) Punched\\Open\\Poisson}{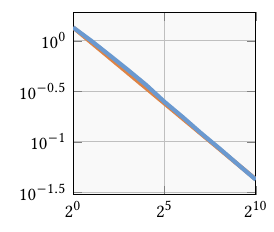}%
    \addgraph{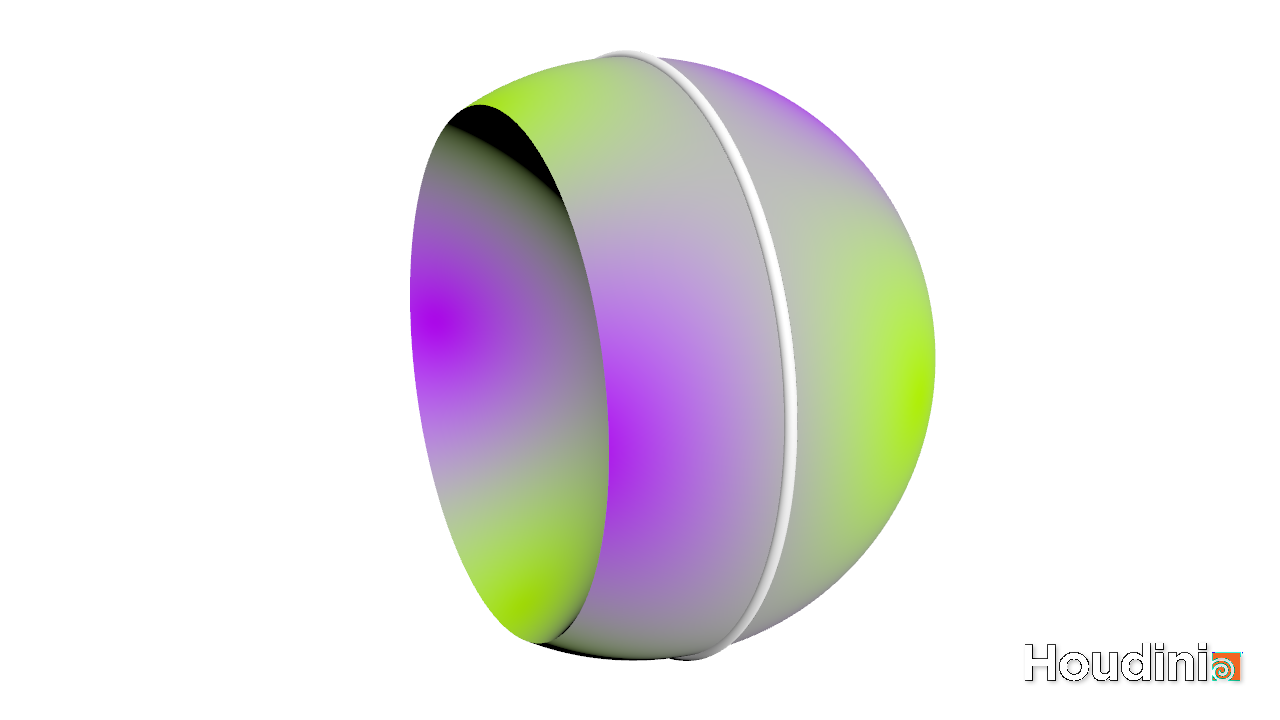}{(o) Punched\\Closed\\Screened Poisson}{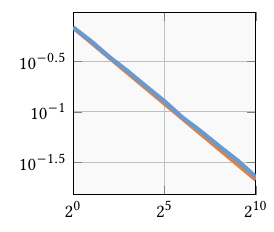}%
    \addgraph{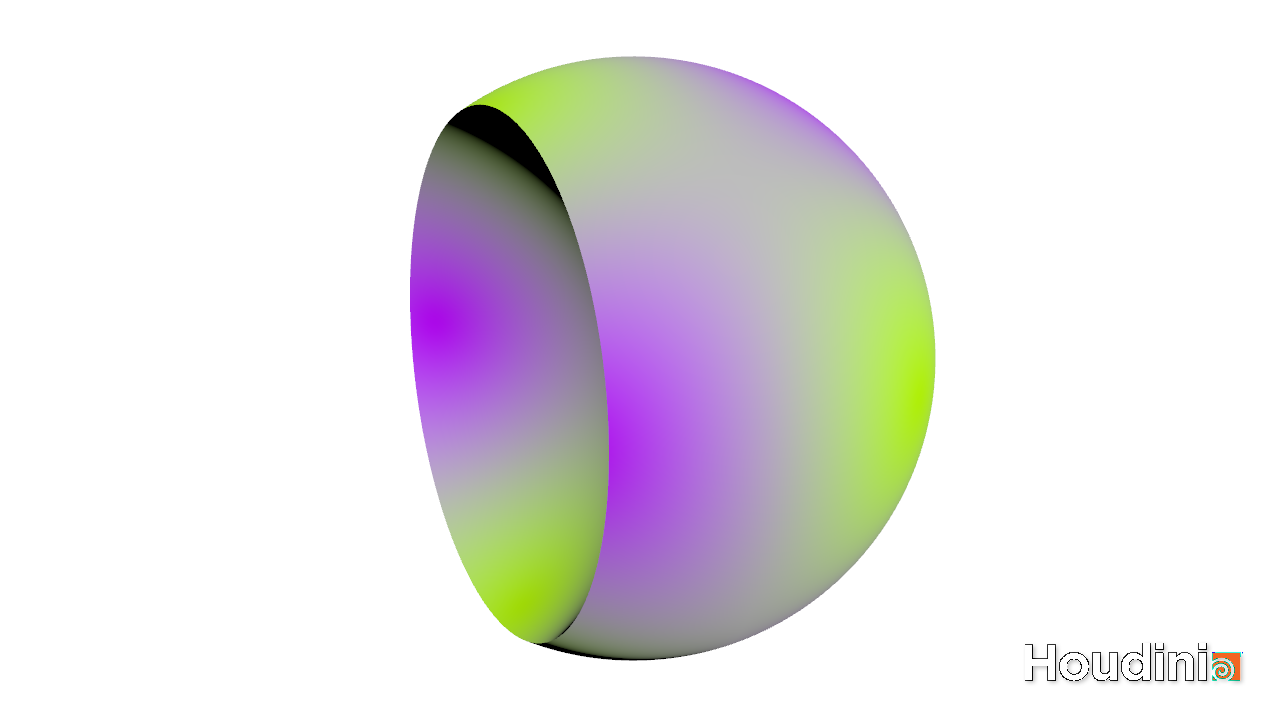}{(p) Punched\\No Boundary\\Screened Poisson}{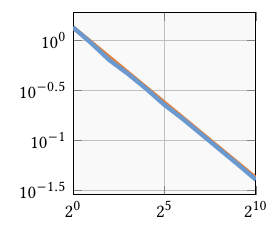}%
    
\caption{Error convergence. In each of the examples, we show the visualization of the scene setup (top), error convergence plot (middle), and scene description (bottom). The scene visualizations show the analytical or reference solution of the problem by mapping the range of solution values to the green-to-purple color gradient and placing a white point or curve on the Dirichlet boundary. The vertical axis of the error plot shows the root mean squared error of the estimates at a few points on the surface in a logarithmic scale, and the horizontal axis shows the number of samples $N_P$ in a logarithmic scale. The blue curves show the result of the experiment, and the orange curves show a line that corresponds to the desired $\mathcal{O}(1/\sqrt{N_P})$ convergence rate. The scene description indicates the surface shape (top), boundary shape/type (middle), and problem type (bottom). While the expected  $\mathcal{O}(1/\sqrt{N_P})$ is achieved in most scenes, the bias remains high when an aggressive medial axis pruning parameter is used (c) and when the source term of the problem is relatively complex (e, g, and h). We describe the details of the scene configurations in the supplemental note.}
\label{fig:convergence}
\end{figure*}

\begin{figure}[H]
\begin{minipage}[ht]{0.40\linewidth}
\centering
\includegraphics[trim={10.5cm 2cm 10.5cm 2cm}, clip, width=0.5\linewidth]{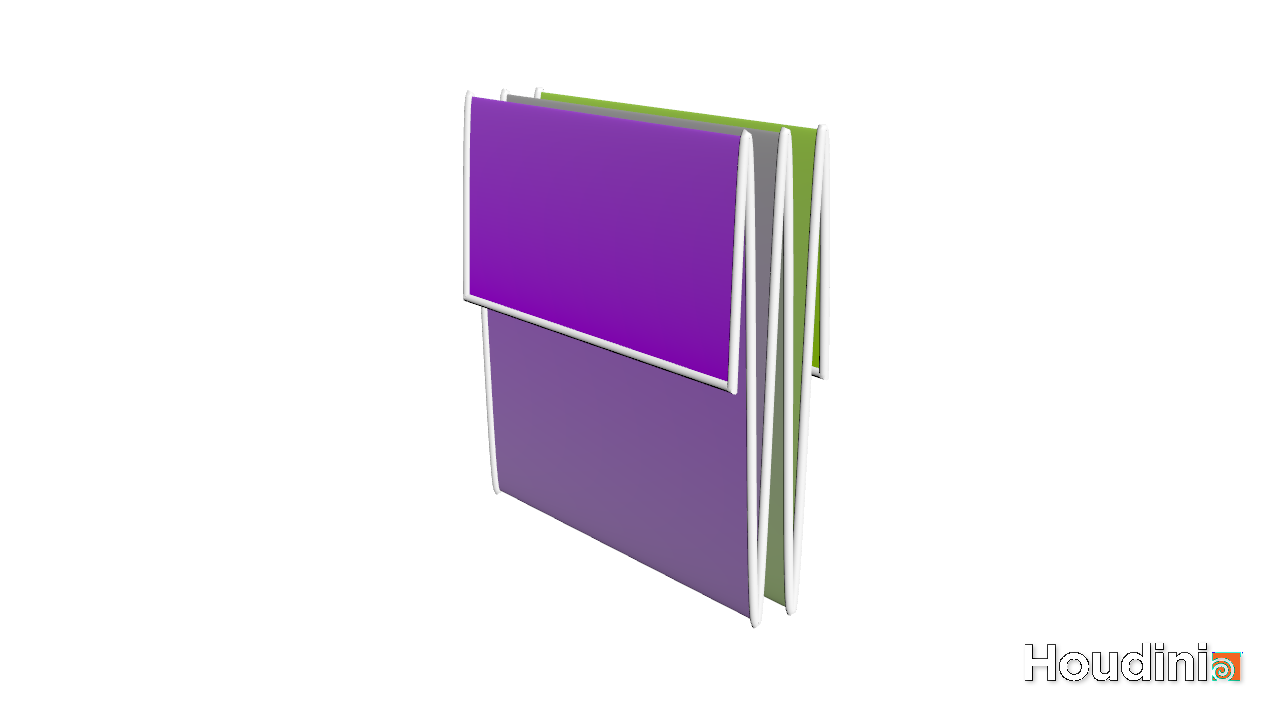}%
\includegraphics[trim={10.5cm 2cm 10.5cm 2cm}, clip, width=0.5\linewidth]{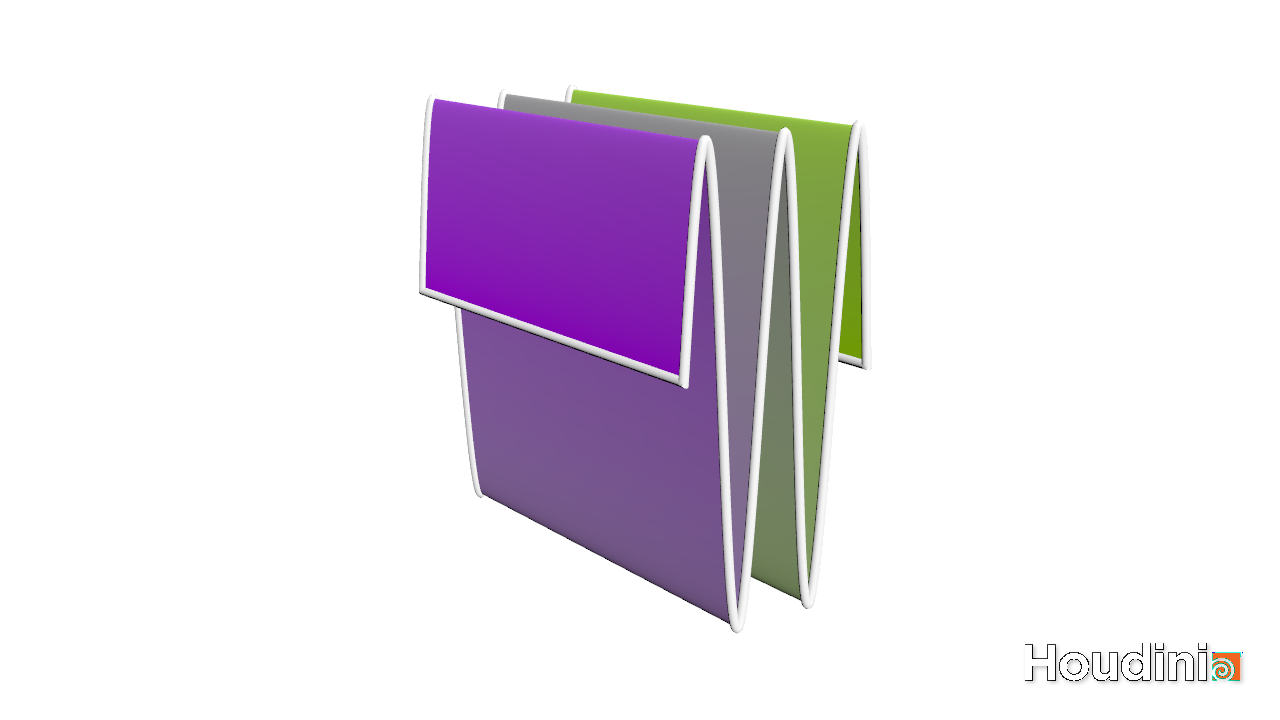}\vspace{-0.5cm}\\
\begin{minipage}[ht]{0.5\linewidth}
\small\textsf{\color{mycolor1}High}
\end{minipage}%
\begin{minipage}[ht]{0.5\linewidth}
\small\textsf{\color{mycolor2}Mid}
\end{minipage}\\
\includegraphics[trim={10.5cm 2cm 10.5cm 2cm}, clip, width=0.5\linewidth]{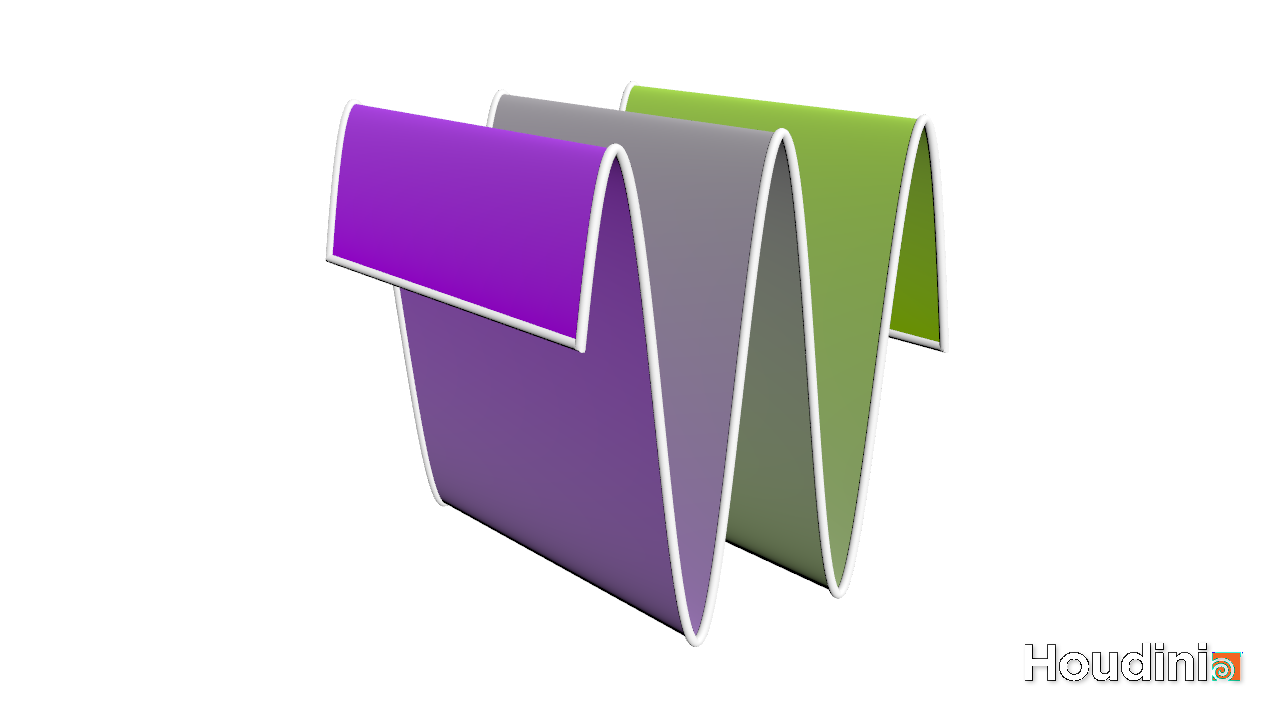}\vspace{-.5cm}\\
\begin{minipage}[ht]{0.5\linewidth}
\small\textsf{\color{mycolor3}Low}
\end{minipage}
\end{minipage}
\hspace{-1cm}
\begin{minipage}[ht]{0.5\linewidth}
\centering
    \resizebox{\linewidth}{!}{\includegraphics{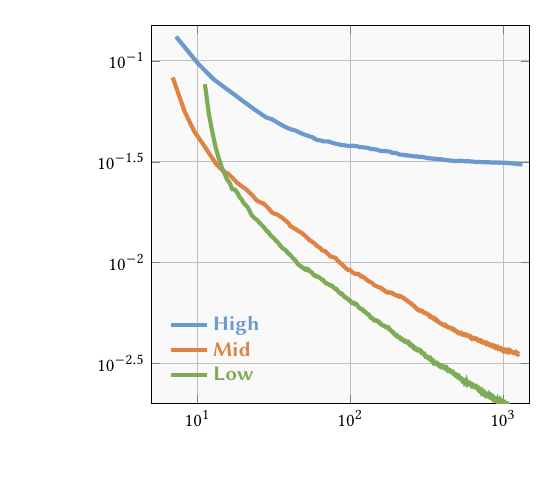}}\\
\end{minipage}
\caption{Effect of local feature size on convergence speed.
We bend a rectangular strip of size 10 by 2 units along sinusoidal curves with high, middle, and low frequencies (top left, top right, and bottom, respectively, in the visualization) and solve the Laplace equation. The analytical solution is defined as the distance along the longer edge of the strip from one of the shorter edges.
The vertical axis of the convergence plot represents the root mean squared error (RMSE), while the horizontal axis shows the time in seconds measured on a MacBook Pro with an M1 Pro chip. We used 1000 sample evaluation points.
The geometry with a larger local feature size allows for faster convergence with lower bias.
}
\label{fig:bentstrip}
\end{figure}

\begin{figure}[H]
\centering
\begin{minipage}[ht]{0.15\linewidth}
\includegraphics[trim={13cm 0cm 19.5cm 0cm}, clip, width=\linewidth]
{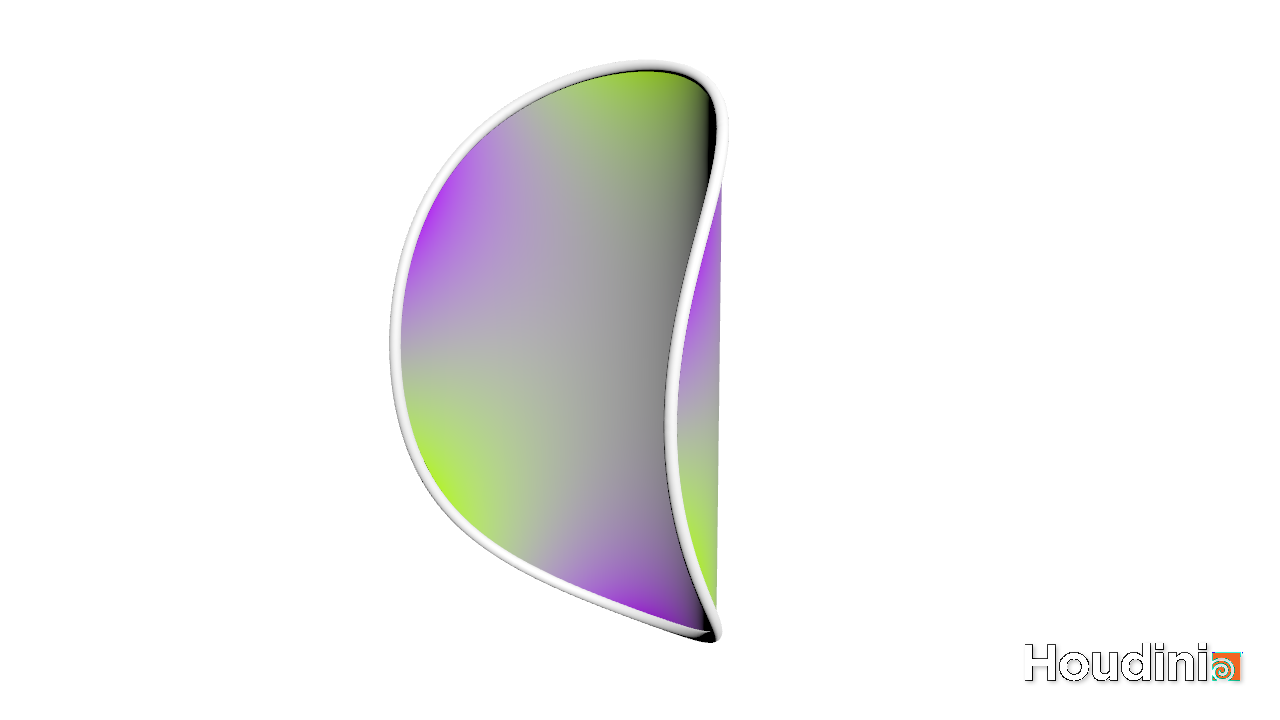}
\end{minipage}%
\begin{minipage}[ht]{0.25\linewidth}
\centering
    \resizebox{1.2\linewidth}{!}{\includegraphics{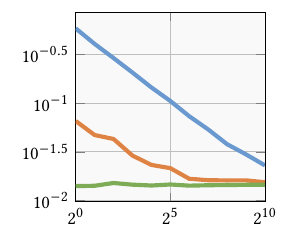}}\\
\end{minipage}
\hspace{0.2cm}
\begin{minipage}[ht]{0.25\linewidth}
\centering
    \resizebox{1.2\linewidth}{!}{\includegraphics{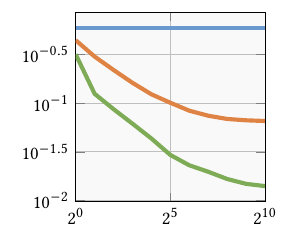}}\\
\end{minipage}
\hspace{0.4cm}
\begin{minipage}[ht]{0.19\linewidth}
\raggedright
\textsf{\textcolor{mycolor1}{-\small No mean \;\;value filtering}}\\
\textsf{\textcolor{mycolor2}{-\small 1x mean \;\;value filtering}}\\
\textsf{\textcolor{mycolor3}{-\small 10x mean \;\;value filtering}}
\end{minipage}\\
\begin{minipage}[ht]{0.17\linewidth}\;\end{minipage}
\begin{minipage}[ht]{0.29\linewidth}\centering\quad\textsf{\small Initial Samples}
\end{minipage}
\begin{minipage}[ht]{0.29\linewidth}\centering\quad\textsf{\small Filter Samples}\end{minipage}
\begin{minipage}[ht]{0.2\linewidth}\;\end{minipage}
\caption{Mean value filtering error. For the Laplace equation with solution $\sol=r^3\sin(3\theta)$ in polar coordinates defined on a unit disk, we curve the disk and measure the error of PWoS with mean value filtering. In each of the two plots, we compare the error without mean value filtering against one application of filtering and ten applications of filtering. The left plot shows the root mean squared error when changing the number of sample paths used to generate the initial estimate with a fixed number of filter samples (1024). The right plot shows the root mean squared error when changing the number of samples used to construct the filter and fixing the number of sample paths used to generate the initial estimate to 1. 
We observe that applying the mean value filter constructed with a sufficiently large number of samples can reduce the error significantly, even when the initial estimate is constructed with a small sample path count.}
\label{fig:convergence_meanvaluefiltering}
\end{figure}

\subsection{Error Convergence}\label{sec:convergence}
We conducted error convergence studies of our method using discretized surfaces. In each scene in \cref{fig:convergence}, we change the number of sample paths $N_P$ to plot the root mean squared error measured against the analytical or reference solution. We do not apply the mean value filtering method in these studies. The scene setup includes Laplace, Poisson, and screened Poisson equations, and both smooth surfaces and surfaces with sharp corners. In all scenes, we use $\epsilon = 0.001$ and $N_V=32$ except for Laplace equations, for which we use $N_V=0$.  While the method shows the expected convergence rate of $\mathcal{O}(1/\sqrt{N_P})$ in most examples, we observed a relatively large bias with problems having more complex source term functions, indicating the need for future work to investigate estimating $g(\vecx)$ from \cref{eq:embedding_pde} for such problems.

\cref{fig:bentstrip} compares the convergence of our method for problems with different local feature sizes. We observe that larger local feature size
corresponds to faster convergence with lower bias.

\paragraph{Mean value filtering.}
We run the mean value filtering algorithm on a Laplace equation on a triangulated curved disk surface in \cref{fig:convergence_meanvaluefiltering}.
In this setup, we observe that applying the filter multiple times with a filter constructed with a sufficiently high sample count can significantly reduce the error, even when the initial estimate is computed with a small number of samples. As expected, the filter is more effective when constructed with more samples, but the error decreases slower than the rate $\mathcal{O}(1/\sqrt{N_P})$, where $N_P$ is the number of samples used to construct the filter. As no recursive estimation is required with the mean value filtering, it significantly improves the efficiency of PWoS.

\begin{figure}[b]
\includegraphics[trim={10.5cm 1.2cm 10.5cm 1.4cm}, clip, width=0.4\linewidth]{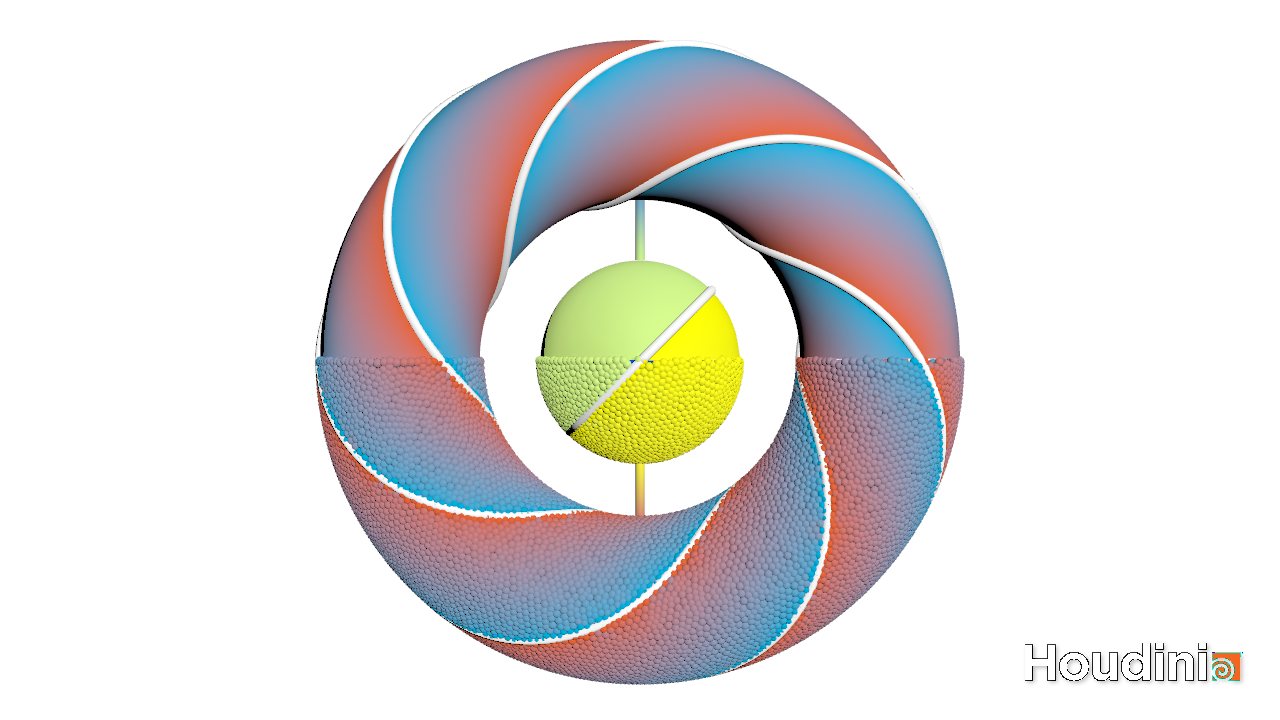}%
\hspace{0.5cm}
\includegraphics[trim={9.5cm 1cm 10cm 1cm}, clip, width=0.4\linewidth]{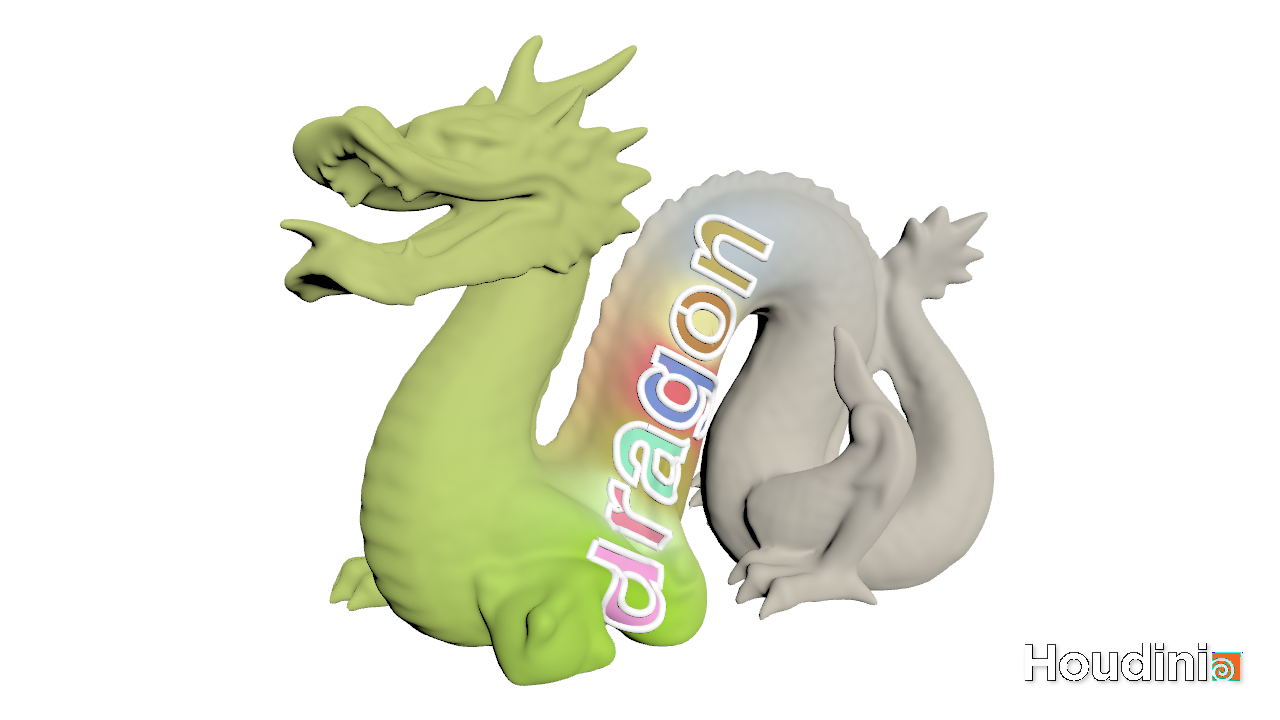}
\caption{Surface diffusion curves solved on various surface representations.
The surface on the left is represented as a combination of triangles, polylines, and oriented points; the surface on the right is represented as a quadrilateral mesh.
The scene on the left, featuring surface geometry of mixed codimension, was adapted from the work of \citet{King:2023:CPM}.}
\label{fig:diffusioncurves}
\end{figure}

\subsection{Applications}\label{sec:applications}

\subsubsection{Diffusion curves}
Diffusion curves \cite{Orzan2008} succinctly represent an image as a collection of curves with associated colors. The final image, exhibiting smooth color gradients, is recovered by solving a Laplace equation with the curves dictating boundary conditions. In our application, we solve the surface Laplace equation using PWoS. 
With our approach, the surface need not have a boundary curve conforming to the discretized mesh, which contrasts with the common approach~\cite{DeGoes2022}.
\cref{fig:diffusioncurves} shows the reconstruction of color at each point on the discretized surface, represented as a quadrilateral mesh and a combination of triangles, polylines, and point clouds.
Our method naturally supports two-sided boundaries, with different colors specified on each side of a curve, and surface geometries with mixed-codimension.
Additionally, the pointwise nature of PWoS allows it to be applied in a view-dependent manner. For example, given a camera configuration, for antialiasing, we sample points within each pixel to generate rays. We then generate PWoS samples at the ray-surface intersection points. No computational resources are wasted on surface points that are invisible to the camera~(\cref{fig:diffusioncurves_viewdependent}), and we can obtain clean results without relying on a fine discretization of the surface.
Boundary integral-based approaches \cite{Bang2023, Sun2012} would similarly allow domain discretization-free evaluation of diffusion curves, but they first require a global linear system solve. Moreover, such methods are not applicable to general curved surfaces, and would need to map the results in the 2D domain to the surface via UV coordinates, for example.

\begin{figure}[b]
\includegraphics[trim={10cm 0.5cm 10cm 0cm}, clip, width=0.25\linewidth]{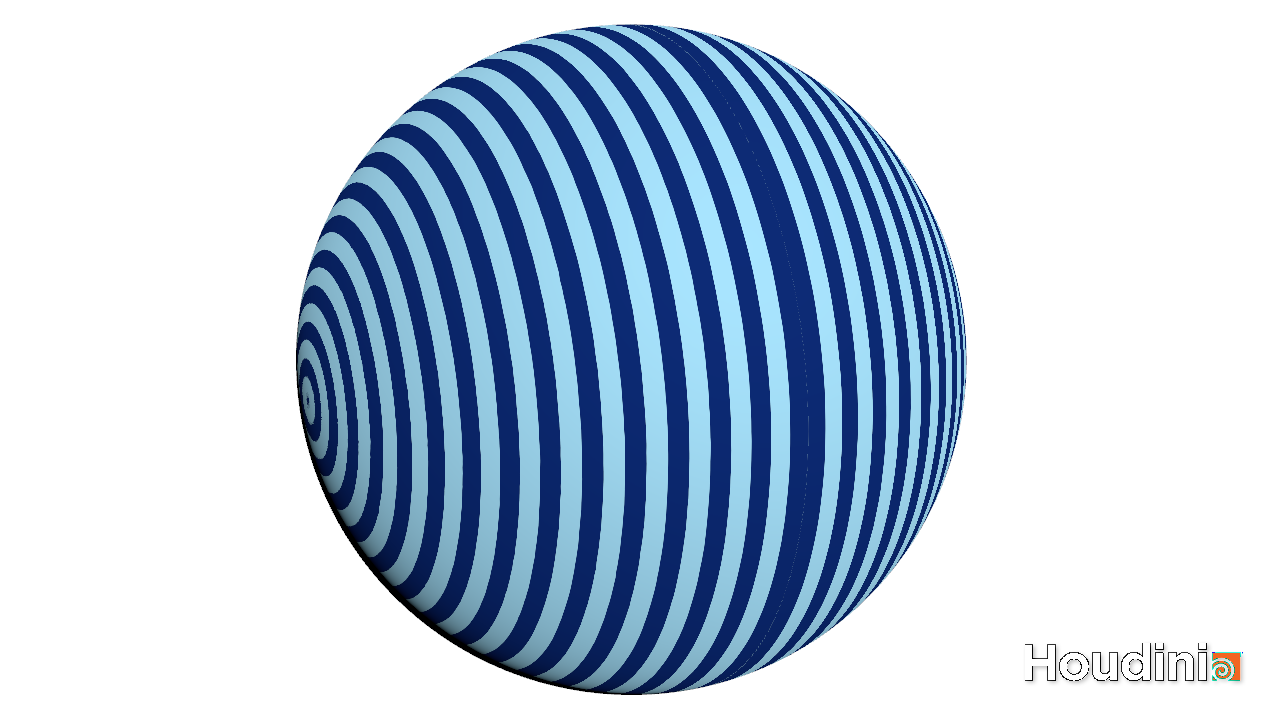}%
\includegraphics[trim={10cm 0.5cm 10cm 0cm}, clip, width=0.25\linewidth]{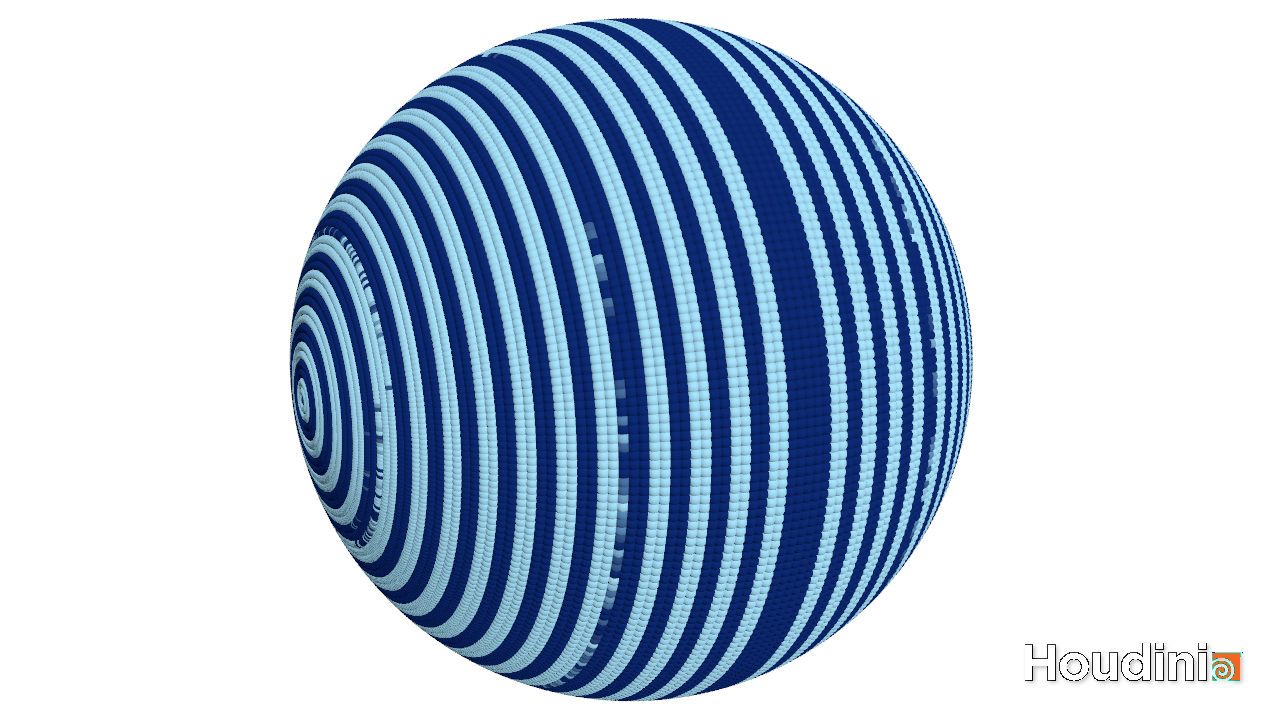}%
\includegraphics[trim={10cm 0.5cm 10cm 0cm}, clip, width=0.25\linewidth]{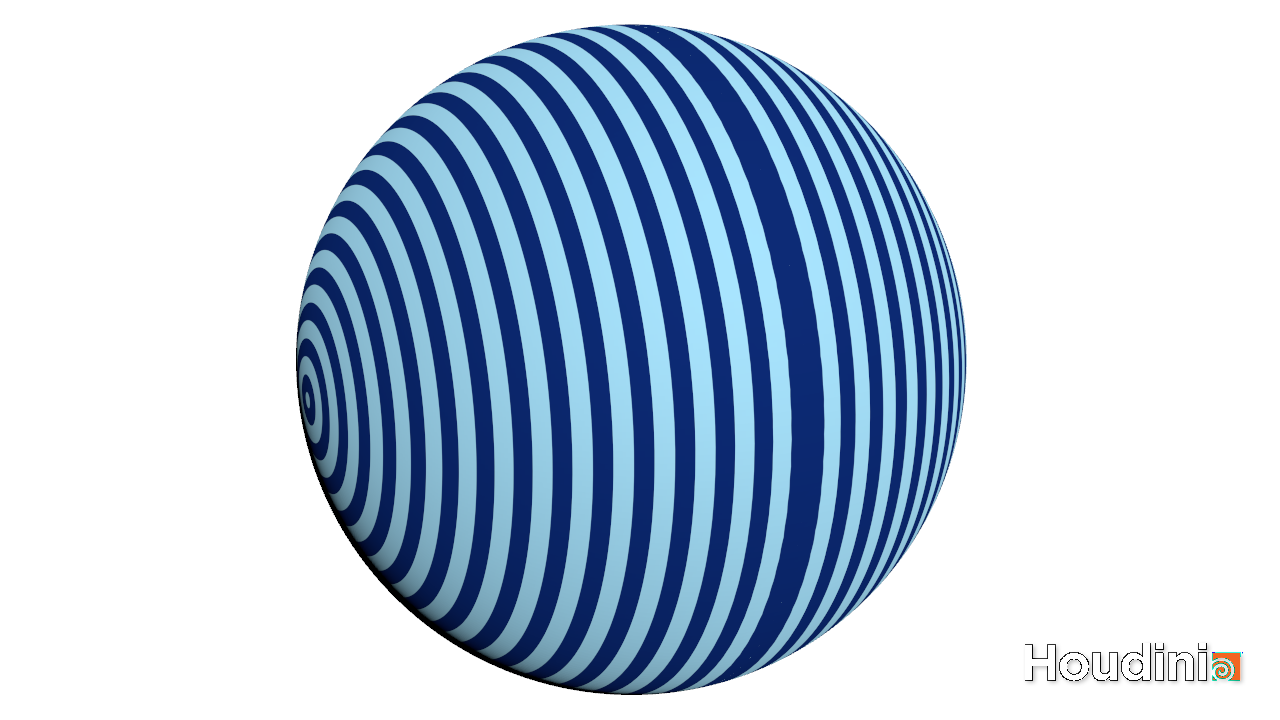}%
\includegraphics[trim={10cm 0.5cm 10cm 0cm}, clip, width=0.25\linewidth]{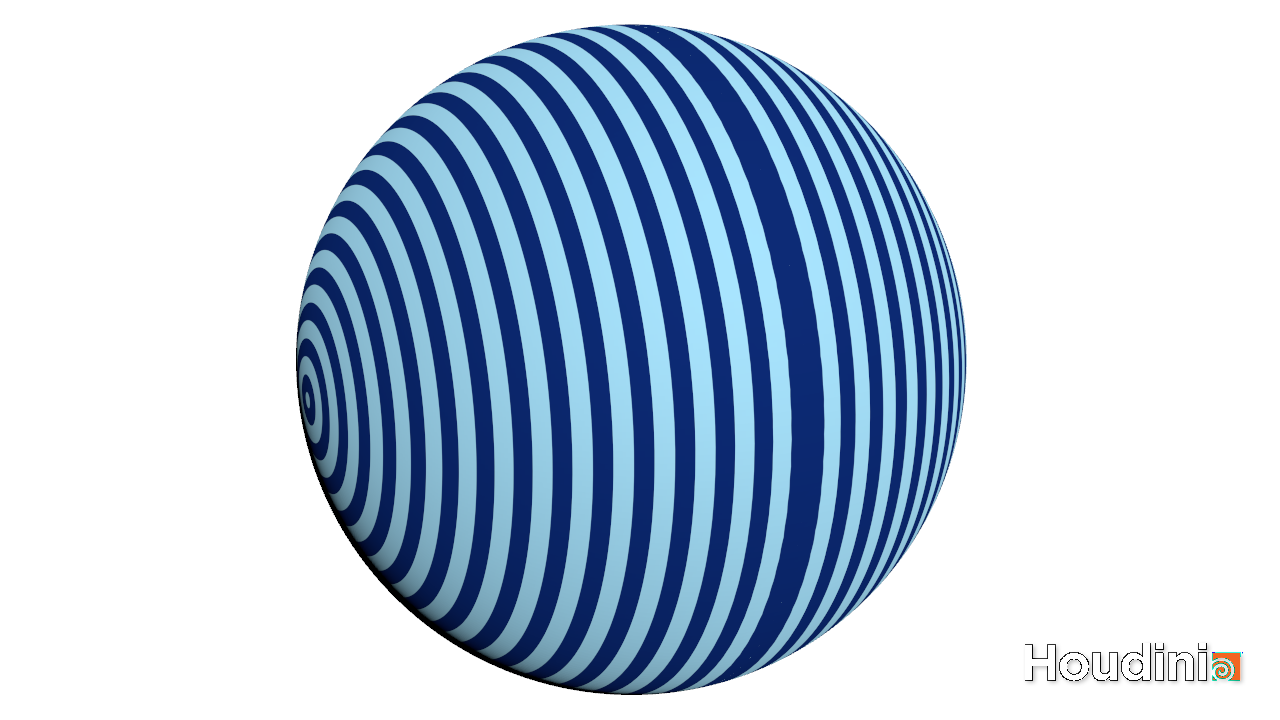}\\
\includegraphics[trim={8cm 1.2cm 12cm 1.4cm}, clip, width=0.25\linewidth]{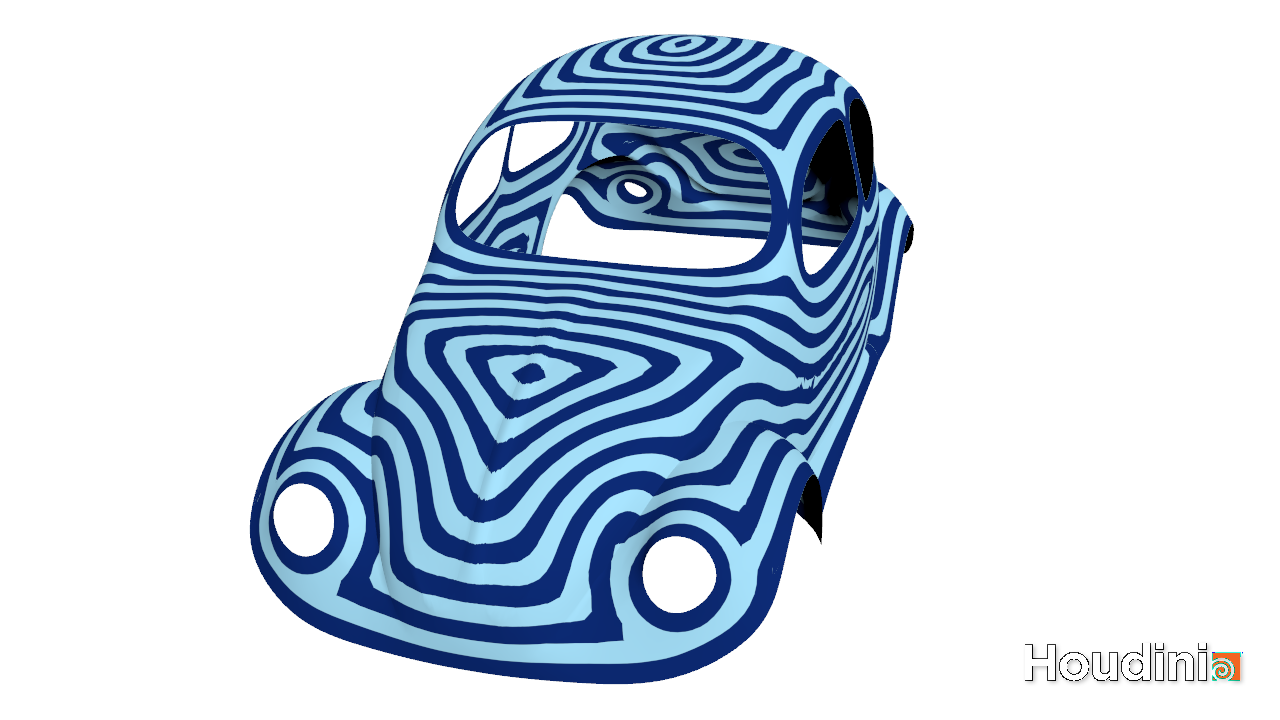}%
\includegraphics[trim={8cm 1.2cm 12cm 1.4cm}, clip, width=0.25\linewidth]{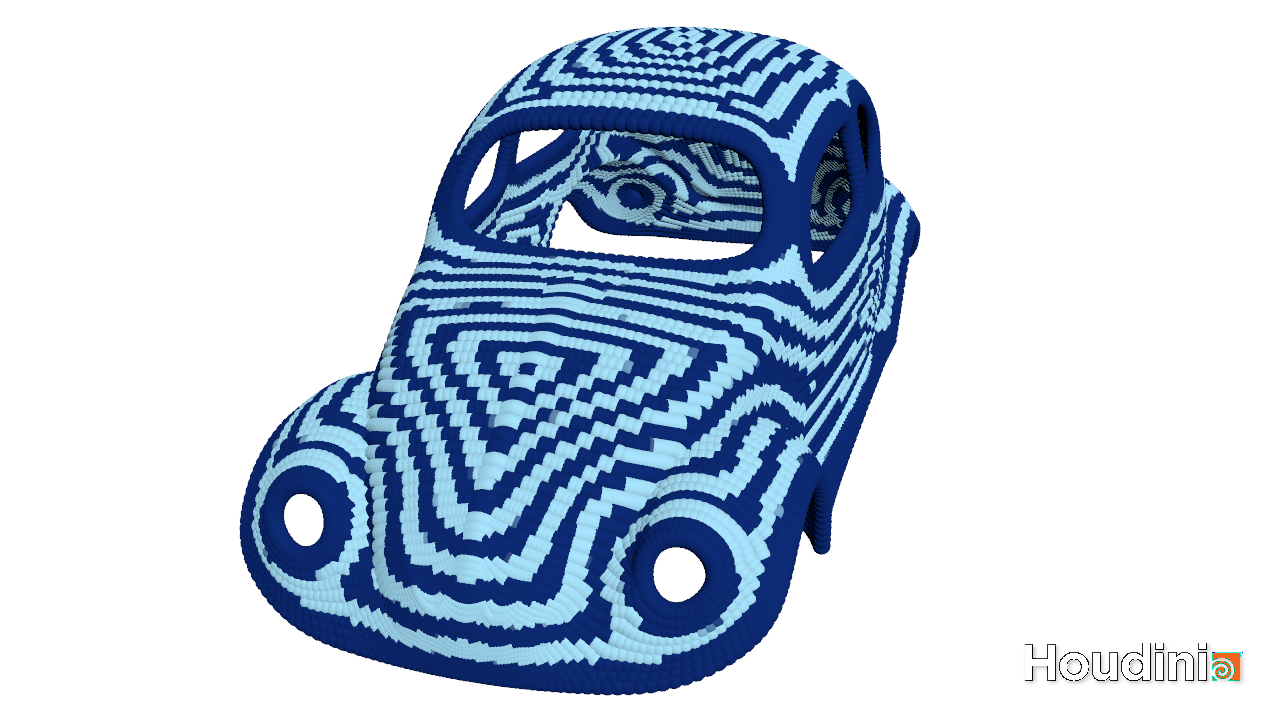}%
\includegraphics[trim={8cm 1.2cm 12cm 1.4cm}, clip, width=0.25\linewidth]{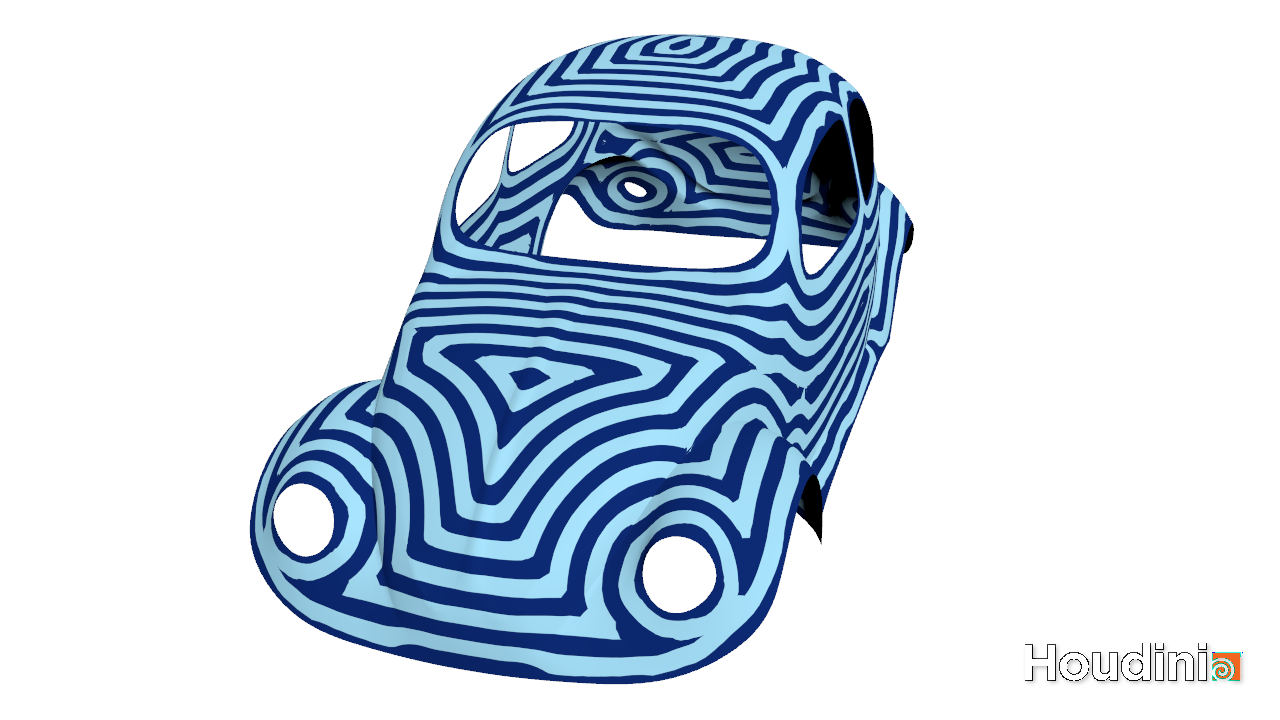}%
\includegraphics[trim={8cm 1.2cm 12cm 1.4cm}, clip, width=0.25\linewidth]{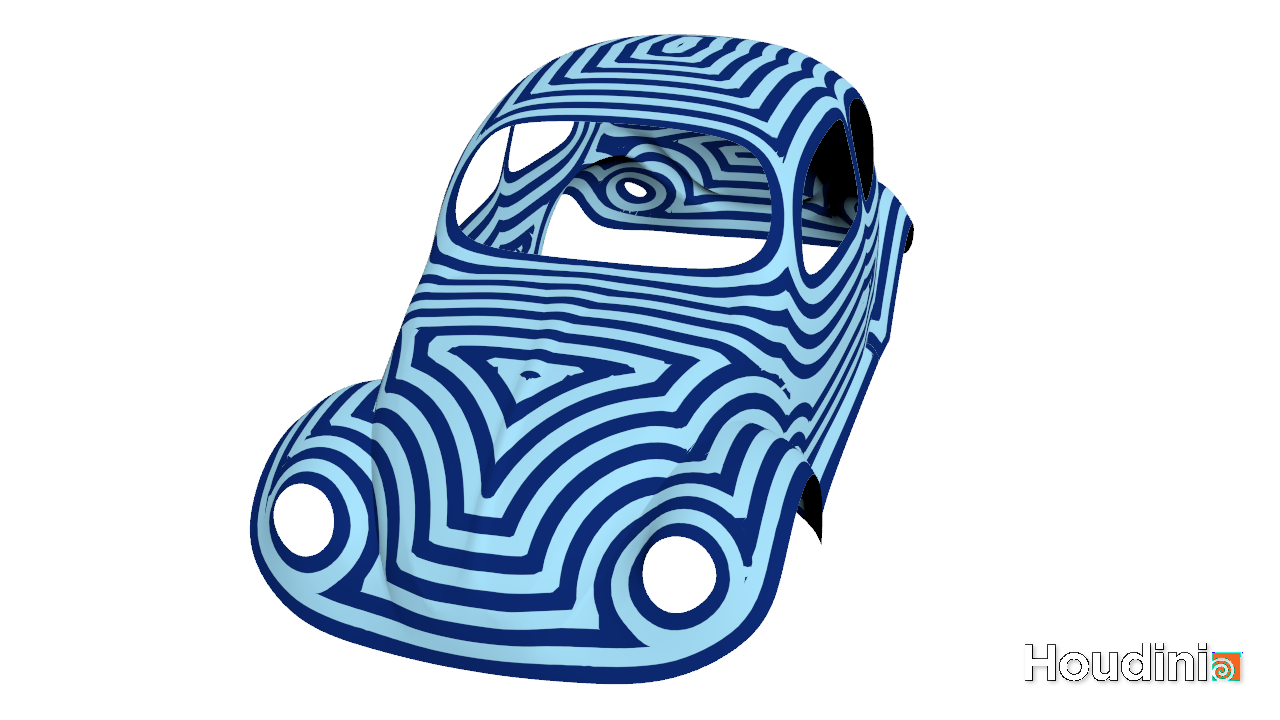}\\
\vspace{0.5em}
\begin{minipage}[hb]{0.25\linewidth}
\centering
    \small\textsf{Triangle Mesh\\(Ours)}
\end{minipage}%
\begin{minipage}[hb]{0.25\linewidth}
\centering
    \small\textsf{Point Cloud\\(Ours)}
\end{minipage}%
\begin{minipage}[hb]{0.25\linewidth}
\centering
    \small\textsf{Grid-Based CPM\\\cite{King:2023:CPM}}
\end{minipage}%
\begin{minipage}[hb]{0.25\linewidth}
\centering
    \small\textsf{Exact\\\cite{geometrycentral}}
\end{minipage}%
\caption{Geodesic distance computation with the heat method. For each of the two scenes, we compare our algorithm on a polygonal mesh representation (leftmost) and oriented point cloud representation (middle-left) against a grid-based CPM counterpart~\cite{King:2023:CPM} (middle-right) and the exact polyhedral distance computed with Geometry Central~\cite{geometrycentral} (rightmost). For the sphere surface (top), we compute the distance from the circle boundary curve in the center, and for the car surface (bottom), we compute the distance from the surface boundary edges. Note that the rendering of the point clouds assigns a UV coordinate per point, resulting in larger visual differences.}
\label{fig:geodesic}
\end{figure}

\subsubsection{Geodesic distance.}\label{sec:geodesic} \citet{Crane2013} proposed the heat method, which solves two standard surface PDEs in series to compute the geodesic distance from the boundary $\Boundary$.
The steps are summarized as
\begin{enumerate}
    \item $\mylaplace_\Surface \sol_\Surface - (1/t)\sol_\Surface = 0$ where $u_\Surface=1\; \mathrm{on}\;\Boundary$,
    \item $\vec{X} = - (\grad_\Surface \sol_\Surface)/\lVert \grad_\Surface \sol_\Surface\rVert_2$, and
    \item $\mylaplace_\Surface \phi_\Surface = \grad_\Surface\cdot \vec{X}$ where $\phi_\Surface=0\; \mathrm{on}\;\Boundary$,
\end{enumerate}
where $t$ is a small positive constant and $\phi$ is the geodesic distance.
Step 2 uses the gradient of the solution to the screened Poisson equation found in Step 1. With a discretization-based method, a discrete gradient operator is used to estimate this gradient; in our method, we directly evaluate the gradient of $\sol$ during Step 1 using the method in \cref{sec:graddiv}, without needing $\sol$ itself. We evaluate the gradient at mesh vertices and normalize it to get $\vec{X}$ at mesh vertices. In Step 3, again, we do not rely on a discrete divergence operator to solve the Poisson equation, but instead use the method described in \cref{sec:graddiv}. When our Poisson solver requires the evaluation of $\vec{X}$ at a point, we interpolate $\vec{X}$ from the mesh vertices and (re)normalize it. We can similarly compute the geodesic distance on a surface represented as a point cloud.
\cref{fig:geodesic} compares our PWoS-based heat method on surfaces represented as polygonal meshes or oriented point clouds against the heat method with grid-based CPM~\cite{King:2023:CPM} and the exact geodesic distance computed with Geometry Central~\cite{geometrycentral}. 
Our results are consistent with the reference implementation, albeit with minor deviations.

\begin{figure}[b]
\includegraphics[trim={10cm 6cm 9cm 3.5cm}, clip, width=0.165\linewidth]{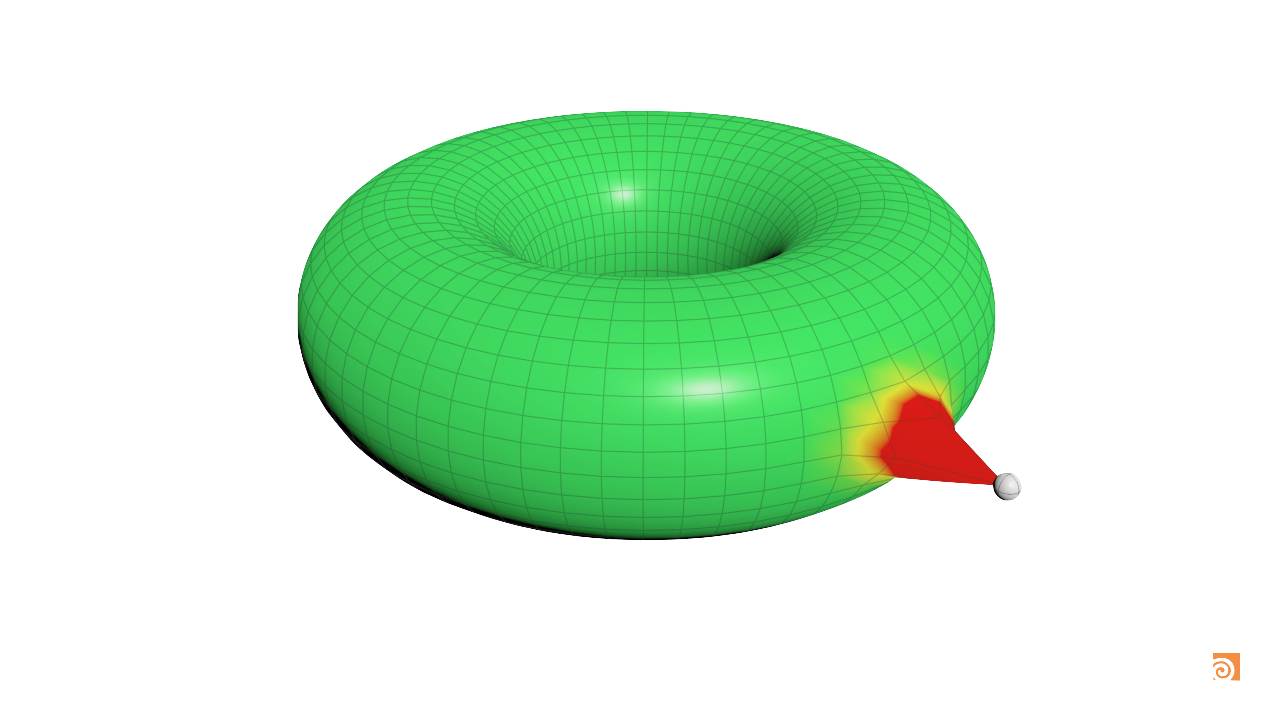}%
\includegraphics[trim={10cm 6cm 9cm 3.5cm}, clip, width=0.165\linewidth]{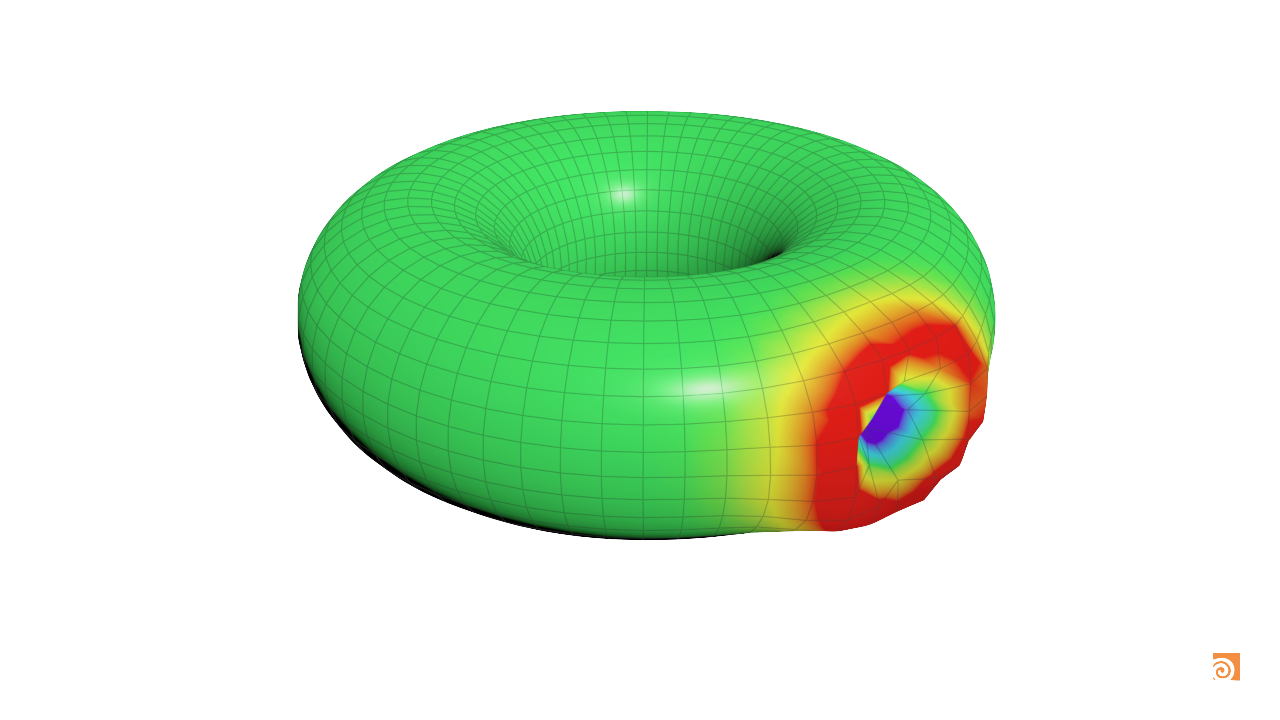}%
\includegraphics[trim={10cm 6cm 9cm 3.5cm}, clip, width=0.165\linewidth]{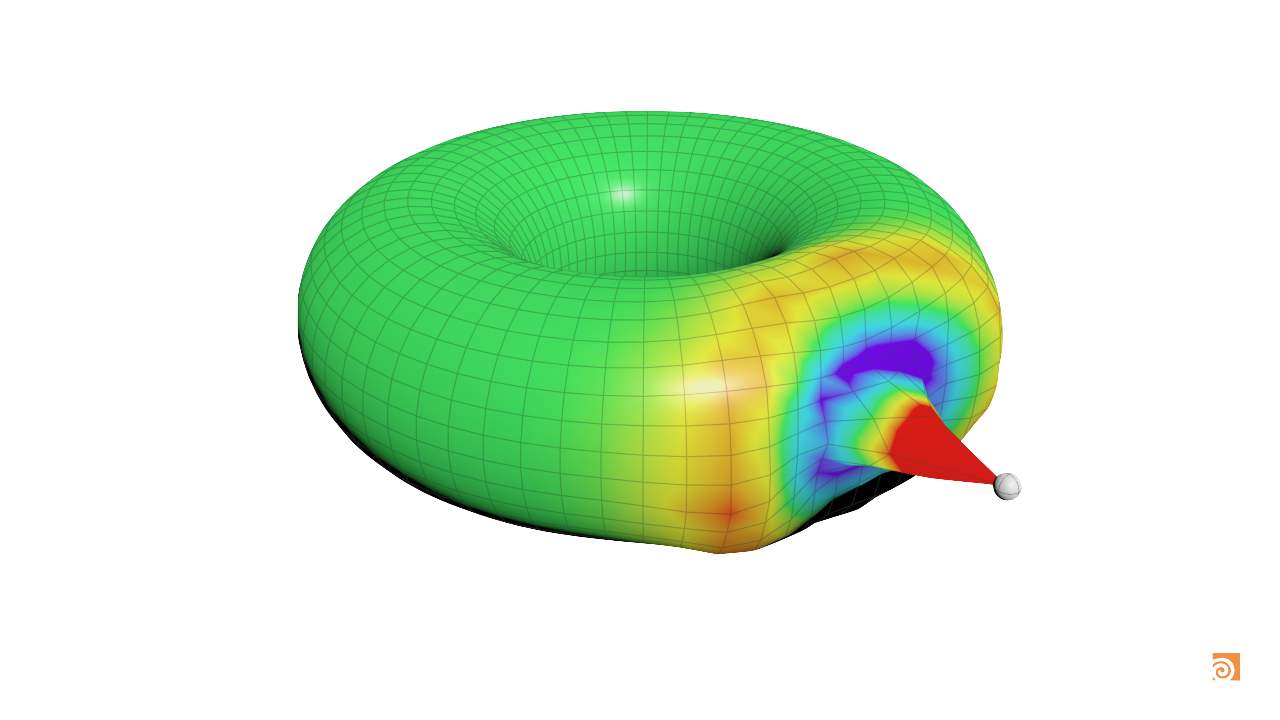}%
\includegraphics[trim={10cm 6cm 9cm 3.5cm}, clip, width=0.165\linewidth]{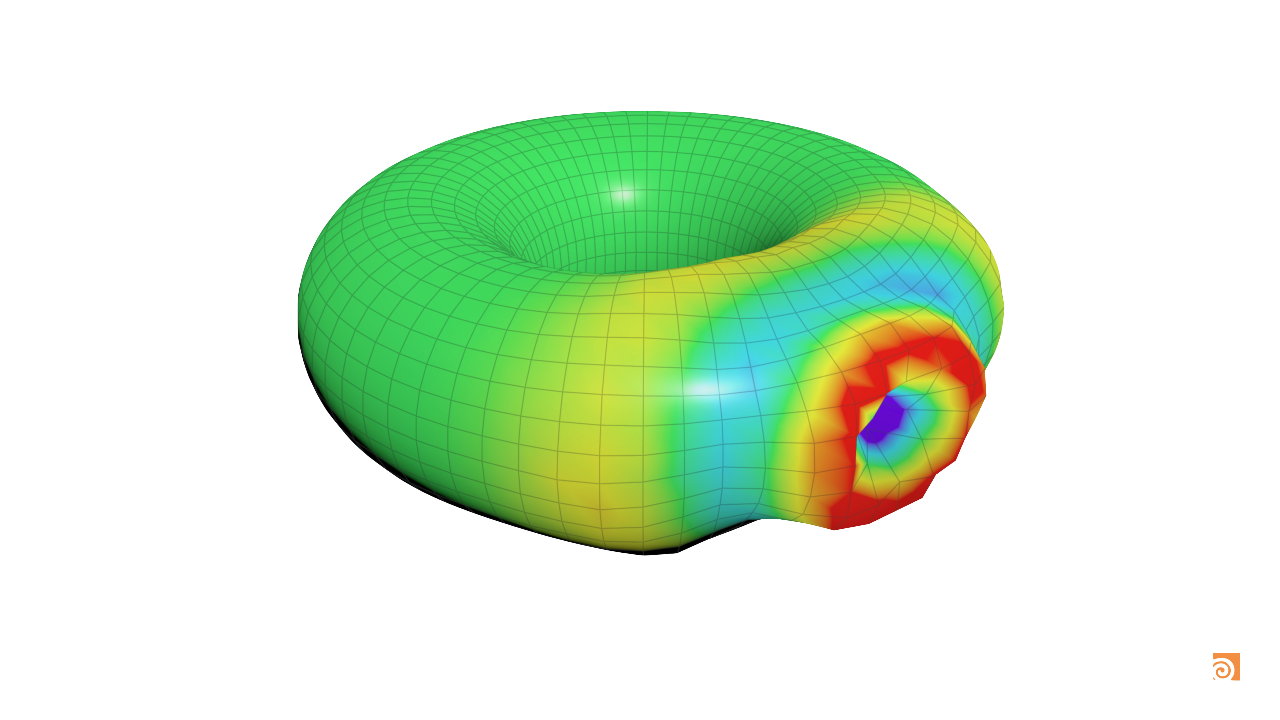}%
\includegraphics[trim={10cm 6cm 9cm 3.5cm}, clip, width=0.165\linewidth]{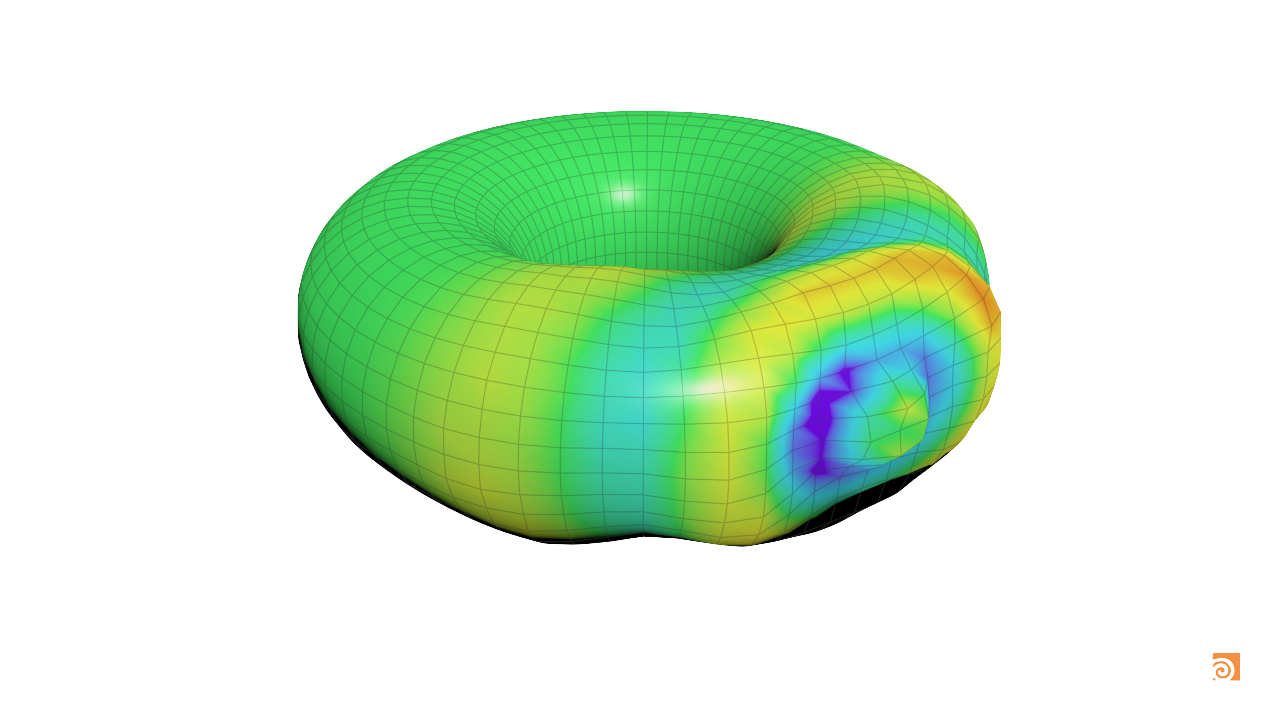}%
\includegraphics[trim={10cm 6cm 9cm 3.5cm}, clip, width=0.165\linewidth]{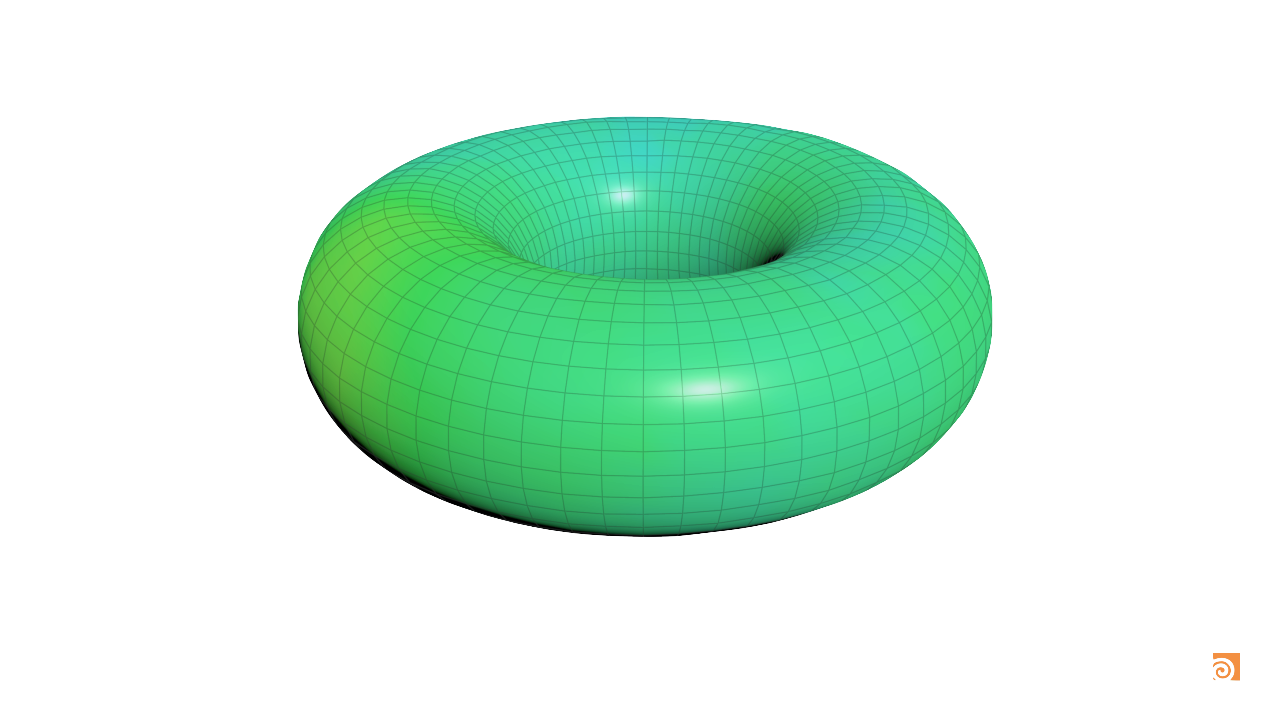}\\
\begin{minipage}[hb]{0.165\linewidth}
\centering
    \small\textsf{Frame 7}
\end{minipage}%
\begin{minipage}[hb]{0.165\linewidth}
\centering
    \small\textsf{Frame 19}
\end{minipage}%
\begin{minipage}[hb]{0.165\linewidth}
\centering
    \small\textsf{Frame 31}
\end{minipage}%
\begin{minipage}[hb]{0.165\linewidth}
\centering
    \small\textsf{Frame 43}
\end{minipage}%
\begin{minipage}[hb]{0.165\linewidth}
\centering
    \small\textsf{Frame 55}
\end{minipage}%
\begin{minipage}[hb]{0.165\linewidth}
\centering
    \small\textsf{Frame 200}
\end{minipage}
\caption{Surface wave animation. We solve the wave equation on the surface and visualize the solution as a displacement applied on the surface in the normal directions. Combined with a time-stepping approach, our method can be applied to a few time-dependent problems. We set the solution to one point on the mesh as the boundary condition for the first 49 frames, and remove the boundary from the 50th frame. The waves damp out as the simulation continues, as expected.}
\label{fig:wave}
\end{figure}

\subsubsection{Surface wave animation}
Using our screened Poisson equation solver, we can solve some classes of time-dependent problems. We discretize the wave equation in time with implicit Euler to get a screened Poisson equation and solve it with time stepping (\cref{fig:wave}). At each time step, we store the solution at the vertices of the mesh and query the solution from previous frames by interpolating the values. In contrast to grid-based CPM~\cite{Auer2012}, our method directly deals with surface geometry without defining an embedding grid.\looseness=-1

\subsection{Performance}
The performance of our method depends on several factors; we report timings for two representative examples. We measured these timings using a workstation with two Intel(R) Xeon(R) Silver 4316 CPUs, each with 20 CPU cores.
For the scene in \cref{fig:diffusioncurves_viewdependent} bottom left, the image resolution is 640 by 480 and the number of samples per pixel was 1024. The precomputation step, including medial axis computation, took less than 1 minute, and the rest of the main parts of PWoS took 2 hours and 11 minutes. We did not apply mean value filtering. 
For the scene in \cref{fig:diffusioncurves} left, we have 28,119 evaluation points. The medial axis point cloud extraction took 2.4 seconds, the initial solution estimate with 1 sample took 1.3 seconds, and the application of 10 mean value filtering steps with a filter constructed with 128 samples took 13 seconds. 
Optimizing the implementation with GPU acceleration may further improve performance. \looseness=-1

\begin{figure}[b]
\hspace{-1cm}
\begin{minipage}[ht]{0.33\linewidth}
\centering
    \resizebox{1.1\linewidth}{!}{\includegraphics{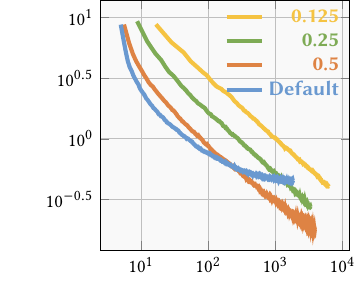}}\\
\end{minipage}%
\begin{minipage}[ht]{0.33\linewidth}
\centering
    \resizebox{1.1\linewidth}{!}{\includegraphics{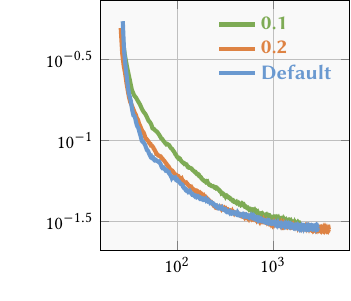}}\\
\end{minipage}%
\begin{minipage}[ht]{0.33\linewidth}
\centering
    \resizebox{1.1\linewidth}{!}{\includegraphics{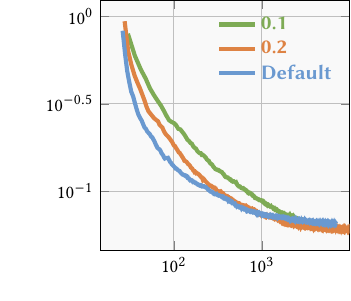}}\\
\end{minipage}\\
\begin{minipage}[ht]{0.33\linewidth}
\centering
\small\textsf{(e)}
\end{minipage}%
\begin{minipage}[ht]{0.33\linewidth}
\centering
\small\textsf{(g)}
\end{minipage}%
\begin{minipage}[ht]{0.33\linewidth}
\centering
\small\textsf{(h)}
\end{minipage}%
\caption{Using bounded sphere size. For the Poisson and screened Poisson problems (e), (g), and (h) in \cref{fig:convergence}, we compare the Default option of not constraining the sphere size (apart from the limit imposed by the local feature size estimate) against specified limits on the maximum sphere size as indicated in the legend. 
The vertical axis shows the root mean squared error, and the horizontal axis shows the time in seconds measured on a MacBook Pro with an M1 Pro chip. For (e), we had 1024 evaluation points, and for (g) and (h), we used 100 sample evaluation points.
While limiting the sphere size may reduce the bias, as we can observe from the intersections of the curves for the default option and the curves for the sphere-size-constrained option, the computation may take longer, and it is difficult to get a practical advantage.
}
\label{fig:lfscap}
\end{figure}

\section{Conclusion and Discussion}
We have developed a Monte Carlo method for surface PDEs by augmenting the formulation of the walk on spheres method with a closest point projection step. Our algorithm is justified through its connection to the theory of closest point extensions drawn from the CPM literature. To accelerate its convergence, we have developed a practical mean value filtering method that utilizes a discrete basis defined over the surface. We have further analyzed the method's convergence on representative analytical tests and demonstrated its application to graphics problems. 

PWoS currently supports only Dirichlet boundary conditions; efficient Neumann or Robin boundary handling similar to the walk on stars method for volumetric PDEs~[\citetalias{Miller:2024:WRW}~\citeyear{Miller:2024:WRW}; \citetalias{sawhney2023walk}~\citeyear{sawhney2023walk}] would require the availability of a few more queries, such as a ray intersection query against the (extended) boundaries. %

While we used a local feature estimation algorithm to allow walks with larger step sizes, the local feature size estimation itself imposes additional smoothness assumptions on the surface. To respect small-scale local features, the walk can require many iterations to reach a Dirichlet boundary. This effect is partly due to our algorithm (like WoS) being based on an integral equation that holds only locally inside a ball. Revisiting this choice using an integral equation based on a global relationship, such as the one underpinning the walk on \emph{boundary} method~\cite{Sugimoto2023:WoB}, could lead to a more efficient alternative for surface PDEs.

Lastly, our method relied on the assumption that the closest point extension compensation term (i.e., $g(\vecx)$ in~\cref{eq:embedding_pde}) in the embedding PDE is negligible. We empirically showed that the algorithm designed with this assumption works well when the source term has a relatively simple expression, but we do not yet have a complete understanding of when this assumption is strictly valid.
However, since $g(\vecx)$ tends to zero continuously as $\vecx$ approaches the surface, the influence of ignoring this term is expected to decrease as we shrink the embedding space (i.e., shrink the sphere
size). One can always take a smaller sphere size, albeit at a higher computational cost, as we show in \cref{fig:lfscap}.
 Extending our method to consider the effect of the compensation term would further improve the reliability and broaden the applicability of our method.

\begin{acks}
This research was partially funded by NSERC Discovery Grants (RGPIN-2021-02524 \& RGPIN-2020-03918), CFI-JELF (Grant 40132), and a grant from Autodesk. 
The first author was partially funded by the David R. Cheriton Graduate Scholarship. 
The second author was partially funded by the Ontario Graduate Scholarship.
\end{acks}

\bibliographystyle{ACM-Reference-Format}
\bibliography{main}

\end{document}

\maketitle

\section{Green's functions and their derivatives}
We list Green's functions on a ball with radius $R$ in $\Rthree$ and their derivatives for readers' convenience.
As \citet{Sawhney:2020:MCGP} summarized, when $\vecx$ is at the center of the ball, the Green's function for the Poisson equation is
\begin{equation}
    G(\vecx, \vecy) = \frac{1}{4\pi}\frac{r-R}{rR},
\end{equation}
and the green's function for the screened Poisson equation is
\begin{equation}
    G_\sigma(\vecx, \vecy) = \frac{1}{4\pi}\frac{\sinh((r-R) \sqrt{\sigma})}{ r \sinh(R \sqrt{\sigma})},
\end{equation}
where $\vecr = \vecy - \vecx$ and $r = \lVert \vecr\rVert_2$.

The gradients of $G$ and $G_\sigma$ with respect to $\vecx$ when $\vecx$ is at the center of the ball are
\begin{equation}
\grad_\vecx G(\vecx, \vecy) = -\frac{1}{4\pi} \left(\frac{1}{r^3} - \frac{1}{R^3}\right) \vecr,
\end{equation}
and
\begin{equation}
\begin{split}
    \grad_\vecx G_\sigma(\vecx, \vecy) =-\frac{1}{4\pi}  &\left(\frac{\sqrt{\sigma} \cosh((R-r) \sqrt{\sigma})}{r\sinh(R\sqrt{\sigma})}\left(\frac{1}{r} - \frac{1}{R}\right)\right.\\
    &\left.+ \frac{\sinh((R-r)\sqrt{\sigma})}{r \sinh(R\sqrt{\sigma}))}
                                \left(\frac{1}{r^2} + \frac{\sqrt{\sigma}\cosh(R\sqrt{\sigma})}{R \sinh(R\sqrt{\sigma})}\right)\right).
\end{split}
\end{equation}

We additionally derive $\grad_\vecx G(\vecx, \vecy)$ in the general case when $\vecx$ is not at the center of the ball:
\begin{equation}
     \grad_\vecx G(\vecx, \vecy) =-\frac{1}{4\pi} \left(\frac{1}{r^3}\vecr + \frac{Ry}{q^3}\vecq\right),
\end{equation}
where $\vecq =  y \vecx - (R^2 /y)\vecy$, $y=\lVert\vecy\rVert_2$, and $q=\lVert\vecq\rVert_2$.
We also have $\grad_\vecz G(\vecx, \vecz) = \grad_\vecz G(\vecz, \vecx)$ due to the symmetry of $G$. We use this expression for problems with a divergence of a vector field as their source term.

\section{Divergence Source Term}
For the solution estimator, when the source term $\source = \diverg\vech$, the volume term converts to
\begin{equation}
\begin{split}
    &\int_{B_\radius(\vecx)} \source(\vecz) G(\vecx, \vecz) \dVz\\
    &=\int_{B_\radius(\vecx)} (\grad_\vecz\cdot\vech(\vecz))\, G(\vecx, \vecz) \dVz,\\
    &=\int_{\partial B_\radius(\vecx)}\vech(\vecz)\cdot\vecn(\vecz)\, G(\vecx, \vecz) \dAz - \int_{B_\radius(\vecx)} \vech(\vecz) \cdot\grad_\vecz G(\vecx, \vecz) \dVz,\\
    &=- \int_{B_\radius(\vecx)} \vech(\vecz) \cdot\grad_\vecz G(\vecx, \vecz) \dVz,
\end{split}
\end{equation}
and we evaluate the last integral instead, which does not require the explicit evaluation of the divergence of $\vech$.
We generate the samples to estimate the converted volume integral with $p(\vecz)\propto 1/\lVert \vecx- \vecz \rVert_2^2$, so the singularity of $\grad_\vecz G$ cancels out.

\section{Gradient Estimation}
The gradient estimator replaces the integral equation for the first step of recursion with
\begin{equation}\label{eq:integraleq_gradient}
    \grad\sol(\vecx) = \frac{1}{\lvert B_\radius(\vecx) \rvert} \int_{\partial B_\radius(\vecx)} \sol(\vecy)\vecn(\vecy) \dAy + \int_{B_\radius(\vecx)} \source(\vecz) \grad_\vecx G(\vecx, \vecz)\dVz,
\end{equation}
where $\lvert B_\radius(\vecx)\rvert$ is the volume of the ball and $\vecn(\vecy)$ is the outward unit normal of the ball at $\vecy$.
For the screened Poisson equation, we multiply the first term by $c_{\radius, \sigma}$ and replace $\grad G$ in the second term with $\grad G_\sigma$ to get a similar integral equation.
To evaluate the integrals, we uniformly sample a point on the sphere for the first term, and we generate the samples with $p(\vecz_i)\propto 1/\lVert \vecx- \vecz_i \rVert_2^2$ for the second term.
The surface gradient of the solution to a surface PDE does not have a normal component, but the estimated solution may have a nonzero normal component before convergence. Thus, to improve the estimate, we set the normal component(s) of the estimated gradient to zero as a post-processing step.

\section{Convergence Study Setup}
We used the following problems to generate the error convergence plots in Fig. 4. Note that we finely discretized the surfaces we describe below to obtain the data we show in the figure.

\paragraph{(a)}
The helix curve we use has three turns, has a radius of $1$, and the endpoints have a height difference of 2. We solve the Laplace equation defined along the curve length $\phi$ as
\begin{equation}
\begin{aligned}
    \frac{\partial^2\sol_\Surface}{\partial\phi^2} &= 0,\\
    u_\Surface(0) &= 0,\\
    u_\Surface(\psi) &= 1,\\
\end{aligned}
\end{equation}
where the boundary conditions are specified at the two ends of the curve, $\phi=0$ and $\phi=\psi$. The analytical solution is $\sol_\Surface(\phi) = \phi/\psi$.

\paragraph{(b) to (d)}
The problem we solve is defined along the curve length $\phi$ as
\begin{equation}
\begin{aligned}
    \frac{\partial^2\sol_\Surface}{\partial\phi^2} &= 0.02,\\
    u_\Surface(0) &= 0,\\
    u_\Surface(\psi) &= 1,\\
\end{aligned}
\end{equation}
where the boundary conditions are specified at the two ends of the curve, $\phi=0$ and $\phi=\psi$, similar to (a). The analytical solution is $\sol_\Surface(\phi) = 0.01\phi^2 + \frac{1 - 0.01\psi^2}{\psi}\phi$.
The helix curve in (b) is identical to the one in (a).
The z-order curve in (c) and (d) is defined using 8 points, $(\pm 1.0, \pm 1.0, \pm 1.0)$.

\paragraph{(e)} This scene is one of the scenes in the grid-based CPM paper by \citet{King:2023:CPM}. On a unit circle, we have a two-sided Dirichlet boundary. In polar coordinates, the problem we solve in terms of the angle $\theta$ is
\begin{equation}
\begin{aligned}
    \frac{\partial^2\sol_\Surface}{\partial\theta^2} &= -2\cos(\theta - \theta_c),\\
    u_\Surface(\theta_c^-) &= 2,\\
    u_\Surface(\theta_c^+) &= 22,\\
\end{aligned}
\end{equation}
where $\theta_c=1.022\pi$ is the position of the Dirichlet boundary. The analytical solution to this problem is $\sol_\Surface(\theta) = 2\cos(\theta - \theta_c) + \frac{10}{\pi}(\theta - \theta_c)$.

\paragraph{(f)} The surface we used is a torus with a major radius $R= 3$ and a minor radius $r=1$. The Dirichlet boundary curve is a torus knot expressed as a parametric curve
\begin{equation}
    x_1(s) = v(s)\cos(as),\;\; x_2(s) = v(s)\sin(as),\;\; x_3(s) = \sin(bs), 
\end{equation}
where $v(s) = R+\cos(bs), a=3, b=7$, and $s\in [0, 2\pi]$. We solve the Laplace equation on the torus with boundary condition $\sin(s)$ along the curve. We used the grid-based CPM implementation of \citet{King:2023:CPM} with a grid spacing of $0.02$ to generate a reference solution and measured the error of PWoS against it.

\paragraph{(g) and (h)}
The surface we used for these setups is the one given by \citet{Dziuk1988}
and later used in multiple CPM works~\cite{Chen2015:CPM,King:2023:CPM}.
The surface is expressed as $\Surface = \{ \vecx \in \Rthree \lvert (x_1-x_3)^2 + x_2^2 + x_3^2=1\}$. The problem we solve is 
\begin{equation}
    \mylaplace_\Surface \sol_\Surface(\vecx) - \sol_\Surface(\vecx)= -\source_\Surface(\vecx), \quad \vecx \in \Surface,
\end{equation}
where $\source_\Surface$ has an analytical, yet complex, expression we can derive as in the work by \citet{Chen2015:CPM}, so the solution of the problem becomes $\sol_\Surface(\vecx) = x_1x_2$. For (g), we use a unit circle on the $x_1x_2$-plane with the analytical solution specified on it as the boundary value as the Dirichlet boundary.
For (h), we did not use any boundary to show the algorithm's convergence for the screened Poisson equation without any boundaries.

\paragraph{(i) to (p)} The scenes consider the unit sphere with a spherical harmonic function as the analytical solution as is done in the study of mesh Laplacians~\cite{Bunge2023}. The sphere mesh is punched inward at $x_3=0.25$ for (m) to (p) to test the algorithm on a geometry with sharp corners.
Given a spherical harmonic $Y_2^3(\vecx) = \frac{1}{4}\sqrt{\frac{105}{\pi}}(x_1^2-x_2^2)^2 x_3$ with eigenvalue $-12$ as the solution, we solve the Poisson equation
\begin{equation}
    \mylaplace_\Surface \sol_\Surface(\vecx) = -12Y_2^3(\vecx), \quad \vecx\in\Surface,
\end{equation}
for (i), (j), (m), and (n) and the screened Poisson equation
\begin{equation}
    \mylaplace_\Surface \sol_\Surface(\vecx) - \sol_\Surface(\vecx) = -13Y_2^3(\vecx), \quad \vecx\in\Surface,
\end{equation}
for (k), (l), (o), and (p).
For (i), (k), (m), and (o), we use the unit circle on the $x_1x_2$-plane as the Dirichlet boundary, and for (j) and (n), we use the unit semicircle where $x_2>0$ as the Dirichlet boundary. We observe the expected convergence behavior with all of the cases in (i) to (p) and suspect that it has something to do with the fact that the source term is a constant multiple of the solution.

\bibliographystyle{ACM-Reference-Format}
\bibliography{main}